\documentclass[a4paper,10pt]{amsart}
\usepackage{epsfig}
\usepackage{graphicx,color}
\usepackage{amsmath,amsfonts,amssymb,eufrak}
\usepackage{float}
\usepackage{graphics}
\usepackage{latexsym}
\usepackage{float}
\usepackage{verbatim}
\usepackage{xy}
\xyoption{matrix}
\xyoption{arrow}
\input xy
\usepackage{rotating}
\usepackage{caption}

\xyoption{all}         

\newtheorem{lem}{Lemma}[section]
\newtheorem{prop}[lem]{Proposition}
\newtheorem{cor}[lem]{Corollary}
\newtheorem{thm}[lem]{Theorem}

\newcommand{\Ext}{\operatorname{Ext}\nolimits}

\newcommand{\mo}{\operatorname{mod}\nolimits}

\newcommand{\End}{\operatorname{End}\nolimits}

\newcommand{\Dim}{\operatorname{dim}\nolimits}

\newcommand{\op}{\operatorname{op}\nolimits}

\begin{document}
\title{The number of elements in the mutation class of a quiver of
  type $D_n$}
\author{Aslak Bakke Buan}
\author{Hermund Andr\' e Torkildsen}

\begin{abstract}
We show that the number of quivers in the mutation class 
of a quiver of Dynkin type $D_n$ is given by $\sum_{d|n}
\phi(n/d)\binom{2d}{d}/(2n)$ for $n \geq 5$. To obtain this formula,
we give a  correspondence between the quivers in the mutation class
and certain rooted trees.
\end{abstract}

\maketitle

\section*{Introduction}
Quiver mutation is an important ingredient in the definition of cluster algebras \cite{fz1}. It is an operation on quivers, which induces an equivalence relation on the set of quivers. The mutation class $\mathcal{M}$ of a quiver $Q$ consists of all quivers mutation equivalent to $Q$. If $Q$ is a Dynkin quiver, then $\mathcal{M}$ is finite. In \cite{t} an excplicit formula for $|\mathcal{M}|$ is given for Dynkin type $A_n$. Here we give an explicit formula for the number of quivers in the 
mutation class of a quiver of Dynkin type $D_n$.
The formula is given by

\[
d(n) = \left\{ \begin{array}{l l}
  \sum_{d|n} \phi(n/d)\binom{2d}{d}/(2n) & \quad \mbox{if $n \geq 5$}, \\
  6 & \quad \mbox{if $n=4$},
\end{array} \right. 
\]
where $\phi$ is the Euler function. 

The proof for this formula consists of two parts.
The first part shows that the mutation class of type 
$D_n$ is in 1--1 correspondence with the triangulations
(with tagged edges) of a punctured $n$-gon,
up to rotation and inversion of tags. This is a generalization of
the method used in \cite{t} to count the number of elements
in the mutation class of quivers of Dynkin type $A_n$.
Here we are strongly using the ideas in \cite{fst} and \cite{s}.

In the second part we count the number of (equivalence classes of)
triangulations of a punctured $n$-gon, by describing an 
explicit correspondence
to a certain class of rooted trees. A tree in this class is constructed by 
taking a family of full binary trees $T_1, \dots, T_s$ such that
the total number of leaves is $n$, and then adding a node $S$ and 
an edge from this node to the root of $T_i$ for each $i$, such that
$S$ becomes a root (Figure \ref{tree tree} displays all such trees for $n=5$).

When these rooted trees are considered up to rotation at the 
root, they are in 1--1 correspondence with the above mentioned
equivalence classes of triangulations of the punctured $n$-gon.
To count these rooted trees we use a simple adaption of a 
known formula found in \cite{i} and \cite[exercise 7.112 b]{st}. 

We also point out a mutation operation on these rooted trees,
corresponding to the other mutation operations involved (on 
triangulations and on quivers).

Our formula and the bijection to triangulations of the punctured
$n$-gon were presented at the ICRA in Torun, August 2007 \cite{t2}.

After completing our work, we learnt about the paper \cite{glz}. They also generalize the methods in \cite{t} to prove the 
bijection from the mutation class of $D_n$ to triangulations of the punctured $n$-gon. However, their
method of counting triangulations is very different from ours. 
They use the classification of quivers of mutation type 
$D_n$, recently given in \cite{v}. 
The authors of \cite{glz} end up with a very different formula than ours.
In particular, their formula is not explicit, and it seems they get a different output than we get,
e.g. for $n=6$.

We are grateful to Hugh Thomas for several useful discussions and for the idea of making use of binary trees as an alternative to rooted planar trees. We would also like to thank Dagfinn Vatne for useful discussions. 

\section{Quiver mutation}
Let $Q$ be a quiver with no multiple arrows, no loops and no
oriented cycles of length two. Mutation of $Q$ at the vertex $k$ gives a
quiver $Q'$ obtained from $Q$ in the following way.

\begin{enumerate}
\item Add a vertex $k^{*}$.
\item If there is a path $i\rightarrow k \rightarrow j$, then if there
  is an arrow from $j$ to $i$, remove this arrow. If there is no arrow
  from $j$ to $i$, add an arrow from $i$ to $j$.
\item For any vertex $i$ replace all arrows from $i$ to $k$ with
  arrows from $k^{*}$ to $i$, and replace all arrows from $k$ to $i$
  with arrows from $i$ to $k^{*}$.
\item Remove the vertex $k$.
\end{enumerate}

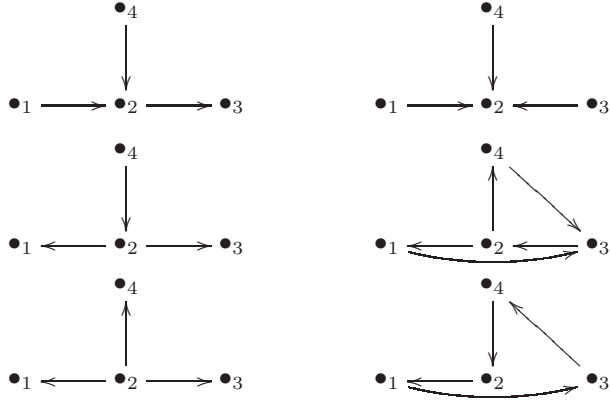
\begin{figure}[htp]
  \begin{center}
	$$\xymatrix{&\bullet_{4}\ar[d]&\\\bullet_{1}\ar[r]&\bullet_{2}\ar[r]&\bullet_{3}} \hspace{1.5 cm}
	\xymatrix{&\bullet_{4}\ar[d]&\\\bullet_{1}\ar[r]&\bullet_{2}&\bullet_{3}\ar[l]}$$
	$$\xymatrix{&\bullet_{4}\ar[d]&\\\bullet_{1}&\bullet_{2}\ar[l]\ar[r]&\bullet_{3}} \hspace{1.5 cm}
	\xymatrix{&\bullet_{4}\ar[dr]&\\\bullet_{1}\ar@/_/[rr]&\bullet_{2}\ar[u]\ar[l]&\bullet_{3}\ar[l]}$$
	$$\xymatrix{&\bullet_{4}&\\\bullet_{1}&\bullet_{2}\ar[r]\ar[l]\ar[u]&\bullet_{3}} \hspace{1.5 cm}
	\xymatrix{&\bullet_{4}\ar[d]&\\\bullet_{1}\ar@/_/[rr]&\bullet_{2}\ar[l]&\bullet_{3}\ar[ul]}$$
  \end{center}\caption{The mutation class of $D_4$.}
  \label{figmutationclassD4} 
  \end{figure}

It is easy to see that mutating $Q$ twice at $k$ gives $Q$. We say
that two quivers $Q$ and $Q'$ are mutation equivalent if $Q'$ can be
obtained from $Q$ by a finite number of mutations. The mutation class
of $Q$ consists of all quivers mutation equivalent to $Q$. Figure
\ref{figmutationclassD4} gives all quivers in the mutation class of
$D_4$, up to isomorphism.  

It is know from \cite{fz3} that the mutation class of a Dynkin quiver
$Q$ is finite. An explicit formula for the number of equivalence
classes in the mutation class of any quiver of type $A_n$ was given in 
\cite{t}. 

The Catalan number $C(i)$ can be defined as the number of
triangulations of an $i+2$-gon with $i-1$ diagonals. It is given by
$$C(i)=\frac{1}{i+1}\binom{2i}{i}.$$ 

The number of equivalence classes in the mutation class of any quiver
of type $A_n$ is then given by the formula \cite{t}

$$a(n)=C(n+1)/(n+3)+C((n+1)/2)/2+(2/3)C(n/3)$$
where the second term is omitted if $(n+1)/2$ is not an integer and
the third term is omitted if $n/3$ is not an integer. This formula counts the triangulations of the disk with $n$ diagonals \cite{b}.

\section{Cluster-tilted algebras}
The cluster category was defined independently in \cite{bmrrt} for the
general case and in \cite{ccs} for the $A_n$ case. Let $\mathcal{D}^b
(\mo H)$ be the bounded derived category of the finitely generated
modules over a finite dimensional hereditary algebra $H$ over a field
$K$. In \cite{bmrrt} the cluster category was defined as the orbit
category $\mathcal{C}=\mathcal{D}^b (\mo H) / \tau^{-1}[1]$, where
$\tau$ is the Auslander-Reiten translation and [1] the suspension
functor. The cluster-tilted algebras are the algebras of the form
$\Gamma=\End_{\mathcal{C}}(T)^{\op}$, where $T$ is a cluster-tilting
object in $\mathcal{C}$ (see \cite{bmr1}). In this paper we will mostly
consider the case where the underlying graph of the quiver of $H$ is
of Dynkin type $D$. 

If $\Gamma = \End_{\mathcal{C}}(T)^{\op}$ for a cluster-tilting object $T$
in $\mathcal{C}$, and $\mathcal{C}$ is the cluster category of a path
algebra of type $D_n$, then we say that $\Gamma$ is of type $D_n$.

Let $Q$ be a quiver of a cluster-tilted algebra $\Gamma$. From \cite{bmr2} it
is known that if $Q'$ is obtained from $Q$ by a finite number of
mutations, then there is a cluster-tilted algebra $\Gamma '$ with quiver
$Q'$. Moreover, $\Gamma$ is of finite representation type if and only if
$\Gamma '$ is of finite representation type \cite{bmr1}. We also have that
$\Gamma$ is of type $D_n$ if and only if $\Gamma '$ is of type
$D_n$. It is well known that we can obtain all orientations of a Dynkin quiver by reflections, and hence all orientations of a Dynkin quiver are mutation equivalent. From \cite{bmr3,birs} we know that a cluster-tilted algebra is up
to isomorphism uniquely determined by its quiver (see also \cite{ccs2}).

It follows from this that the number of non-isomorphic cluster-tilted
algebras of type $D_n$ is equal to the number of equivalence classes in the
mutation class of any quiver with underlying graph $D_n$. 

\section{Category of diagonals of a regular $n+3$-gon}
In \cite{ccs} Caldero, Chapoton and Schiffler considered regular polygons with $n+3$ vertices and triangulations of such polygons. A diagonal is a
straight line between two non-adjacent vertices on the border of the polygon,
and a triangulation is a maximal set of diagonals which do not
cross. A triangulation of an $(n+3)$-gon consists of exactly $n$
diagonals. In \cite{ccs} the category of diagonals of such polygons was
defined, and it was shown to be equivalent to the cluster category, as
defined in Section 2, in the $A_n$ case. It was also shown that a
cluster-tilting object in the cluster category $\mathcal{C}$ corresponds to a
triangulation of the regular $(n+3)$-gon in the $A_n$ case. In
\cite{t} it was shown that there is a bijection between isomorphism
classes of cluster-tilted algebras of type $A_n$ (or equivalently
isomorphism classes of quivers in the mutation class of any quiver with
underlying graph $A_n$) and triangulations of the disk with $n$
diagonals (i.e. triangulations of the regular $(n+3)$-gon up to
rotation).  

For any triangulation of the regular $(n+3)$-gon we can define a
quiver with $n$ vertices in the following way. The vertices are the
midpoints of the diagonals. There is an arrow between $i$ and $j$ if
the corresponding diagonals bound a common triangle. The orientation
is $i \rightarrow j$ if the diagonal corresponding to $j$ can be
obtained from the diagonal corresponding to $i$ by rotating
anticlockwise about their common vertex. It is also known from
\cite{ccs} that all quivers obtained in this way are quivers of
cluster-tilted algebras of type $A_n$. This means that we can define a
function $\gamma_n$ from the mutation class of $A_n$ to the set of all
triangulations of the regular $(n+3)$-gon. There is an induced
function $\widetilde{\gamma}_n$ from the mutation class of $A_n$ to the set
of all triangulations of the disk with $n$ diagonals. It was shown in
\cite{t} that $\widetilde{\gamma}_n$ is a bijection. 
\begin{figure}[htp]
  \begin{center}
    \includegraphics[width=5.0cm]{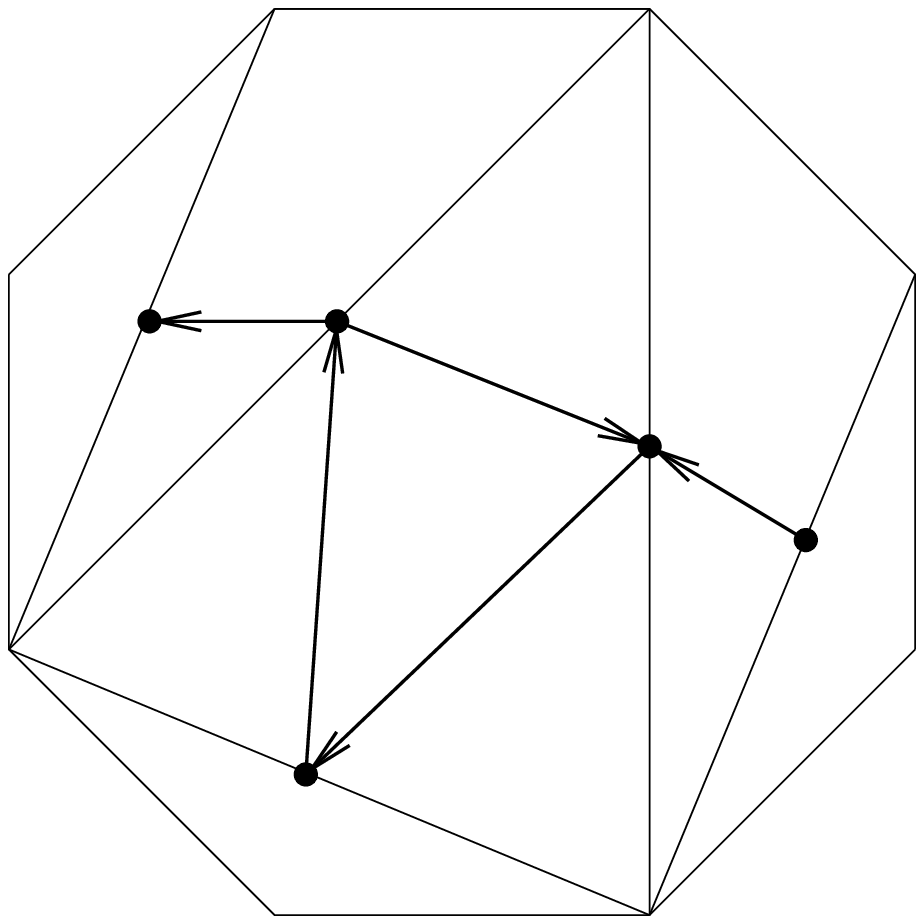}
  \end{center}\caption{A triangulation $\Delta$ of the regular 8-gon and the
    corresponding quiver $\gamma_5(\Delta)$ of type $A_5$.}
  \label{figeksan} 
  \end{figure}

\section{Category of diagonals of a punctured regular $n$-gon}
In this paper we will consider the $D_n$ case and we will first recall some results and notions from \cite{s} and \cite{fst}.

Let $\mathcal{P}_n$ be a regular polygon with $n$ vertices and one puncture in the
center. Diagonals (or edges) will be homotopy classes of paths between
two vertices on the border of the polygon. We follow the definitions
from \cite{s}.

Let $\delta_{a,b}$ be an oriented path between two vertices $a \neq b$ on the
border of $\mathcal{P}_n$ in counterclockwise direction, such that
$\delta_{a,b}$ does not run through the same point twice. Also let
$\delta_{a,a}$ be the path that runs from $a$ to $a$, i.e. around
the polygon exactly one time. We define $|\delta_{a,b}|$ to be the
number of vertices on the path $\delta_{a,b}$, including $a$ and
$b$. 

An edge is a triple $(a,\alpha,b)$ where $a$ and $b$ are vertices
on the border of the polygon and $\alpha$ is an oriented path from $a$
to $b$ lying in the interior of $\mathcal{P}_n$ and that is homotopic
to $\delta_{a,b}$. Furthermore, the path should not cross itself and
$|\delta_{a,b}| \geq 3$. Two edges are equivalent if they start in the
same vertex, end in the same vertex and are homotopic. 

Let $E$ be the set of equivalence classes of edges, and denote by
$M_{a,b}$ the equivalence class of edges in $E$ going from $a$ to
$b$. In \cite{s} the set of tagged edges is defined as follows.  

$$\{ M_{a,b}^{\epsilon}|M_{a,b} \in E\text{, } \epsilon \in \{-1,1\} \text{
  with }
\epsilon = 1 \text{ if } a \neq b \}$$

From now on tagged edges will be called \textit{diagonals}. Diagonals
starting and ending in the same vertex $a$ will be represented as
lines between the puncture and the vertex $a$. Diagonals with
$\epsilon = -1$ will be drawn with a tag on it. In some cases we will
draw them as loops.  

The crossing number $e(M_{a,b}^\epsilon,N_{c,d}^{\epsilon'})$ is the
minimal number of intersection of representations of
$M_{a,b}^\epsilon$ and $N_{c,d}^{\epsilon'}$ in the interior of the
punctured polygon. When $a=b$ and $c=d$, we let the crossing
number be $1$ if $a \neq c$ and $\epsilon \neq \epsilon'$ and $0$
otherwise. If $e(M_{a,b}^\epsilon,N_{c,d}^{\epsilon'}) = 0$, we say
that $M_{a,b}^\epsilon$ and $N_{c,d}^{\epsilon'}$ do not cross. 

\begin{figure}[htp]
  \begin{center}
    \includegraphics[width=4.0cm]{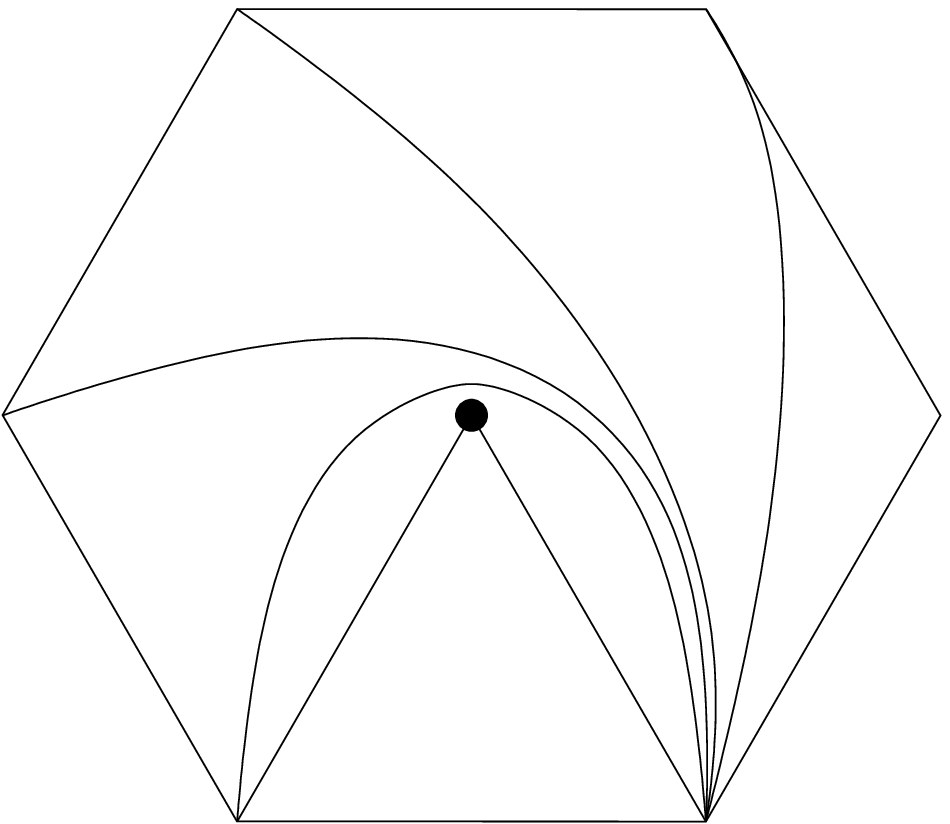}
    \includegraphics[width=4.0cm]{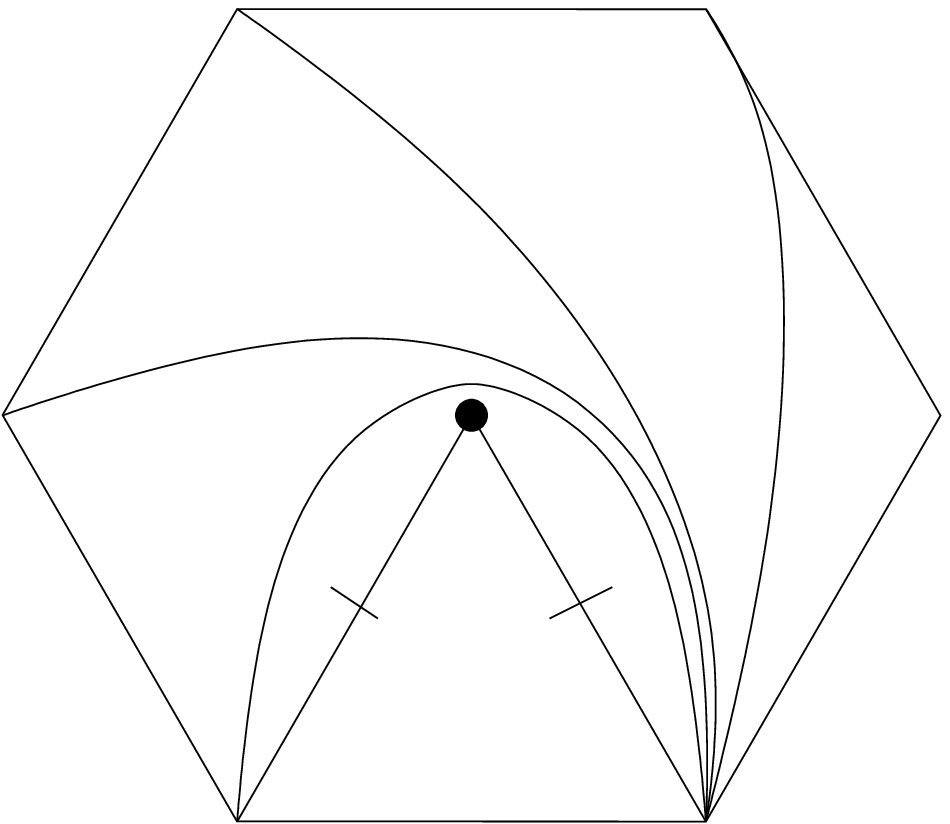}
    \includegraphics[width=4.0cm]{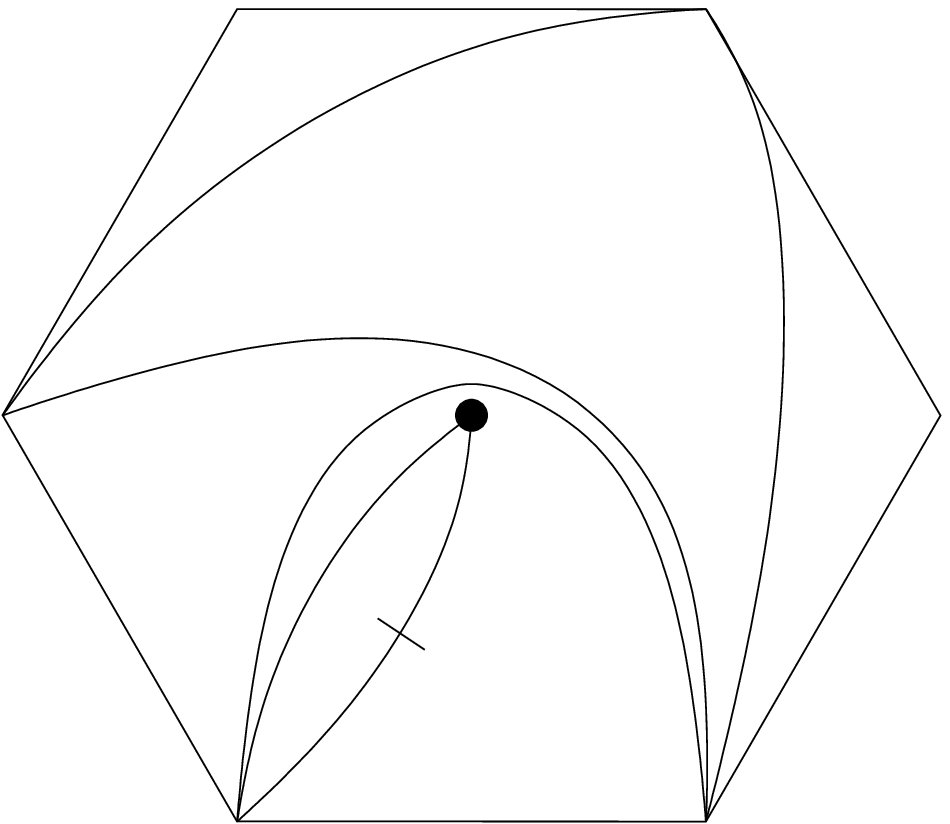}
  \end{center}\caption{Examples of triangulations of the punctured 6-gon.}
  \label{figekstagged} 
  \end{figure}

Now we can define a triangulation of the punctured $n$-gon, which
is a maximal set of non-crossing diagonals. Any such set will have $n$
elements \cite{s}. See some examples of triangulations of the
punctures 6-gon in Figure \ref{figekstagged}.

\cite{s} defines a category which is equivalent to the cluster
catecory in the $D_n$ case in the following way. The objects are
direct sums of diagonals (tagged edges), and the morphism space from
$\alpha$ to $\beta$ is spanned by sequences of elementary moves modulo
the mesh-relations. The equivalence between this category
$\mathcal{C}$ and the cluster category in the $D_n$ case was proved in
\cite{s}. Furthermore we have the following important results: 

\begin{itemize}
\item $\Dim \Ext_{\mathcal{C}}^1(\alpha,\beta)$ is equal to the crossing
  number of $\alpha$ and $\beta$.
\item A cluster-tilting object corresponds to a triangulation.
\item The Auslander-Reiten translation of a diagonal from $a$ to $b$
  is given by clockwise rotation of the diagonal if $a \neq b$. If
  $a=b$ the AR-translation is given by clockwise rotation and
  inverting the tag. 
\end{itemize}

\begin{figure}[htp]
  \begin{center}
    \includegraphics[width=4.0cm]{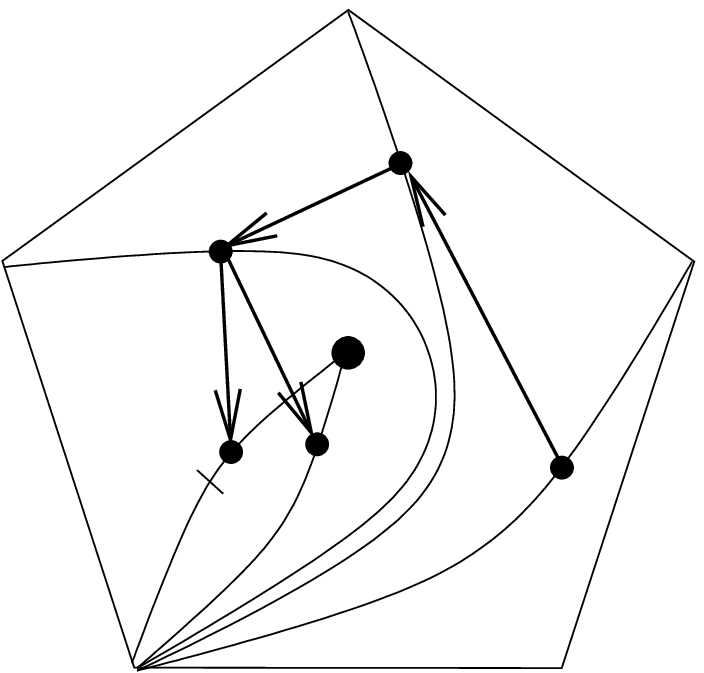}
    \includegraphics[width=4.0cm]{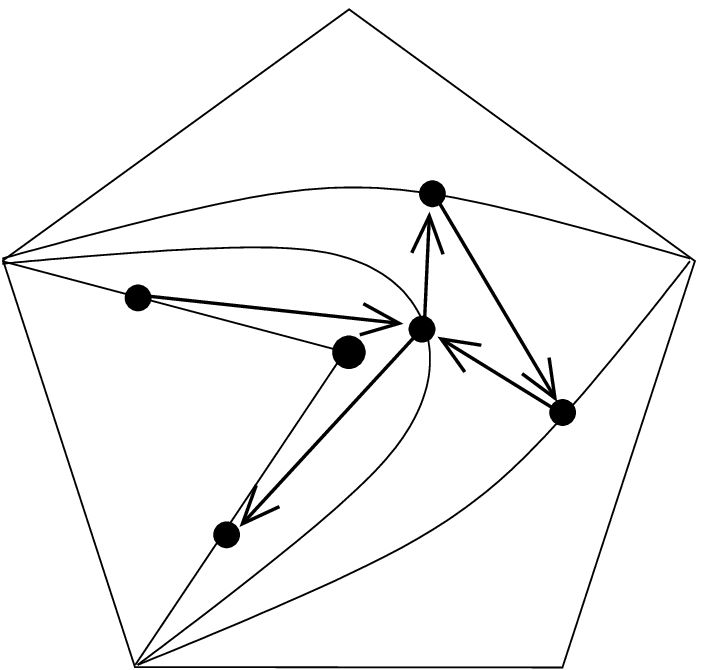}
  \end{center}\caption{Some examples of triangulations and corresponding quiver.}
  \label{figtriquiv} 
  \end{figure}

Let $\mathcal{T}_n$ be the set of all triangulations of
$\mathcal{P}_n$, and let $\Delta$ be an element in $\mathcal{T}_n$. We can assign to $\Delta$ a quiver in the following way (see \cite{fst}). Just as in the $A_n$ case, the vertices are the
midpoints of the diagonals. There is an arrow between $i$ and $j$ if
the corresponding diagonals bound a common triangle. The orientation
is $i \rightarrow j$ if the diagonal corresponding to $j$ can be
obtained from the diagonal corresponding to $i$ by rotating
anticlockwise about their common vertex. In the case when there are two diagonals $\alpha$ and $\alpha'$ between the puncture and the same vertex on the border, both adjacent to a diagonal $\beta$ and a border edge $\delta$, we consider the triangle with edges $\alpha$, $\beta$ and $\delta$ separately from the triangle with edges $\alpha'$, $\beta$ and $\delta$, when thinking of $\alpha$ and $\alpha'$ as loops around the puncture. If we end up with an oriented cycle of length $2$, delete both arrows in the cycle. See some examples in Figure \ref{figtriquiv}.

Let $\mathcal{M}_n$ be the mutation class of $D_n$, i.e. all quivers
obtained by repeated mutations from $D_n$, up to isomorphisms of
quivers. We can define a function $\epsilon_n :
\mathcal{T}_n \rightarrow \mathcal{M}_n$, where we set $\epsilon_n(\Delta) = Q_{\Delta}$ for any triangulation in $\mathcal{T}_n$. It is known that $Q_{\Delta}$ is a quiver of Dynkin type $D_n$ and that all quiver of type $D$ can be obtained this way, hence $\epsilon$ is surjective. 

We can define a mutation operation on a triangulation. If $\alpha$ is a diagonal in a triangulation, then mutation at $\alpha$ is defined as replacing $\alpha$ with another diagonal such that we obtain a new triangulation. This can be done in one and only one way. It is known that mutation of quivers
commutes with mutation of triangulations under $\epsilon$ (see \cite{s,fst}).

\section{Bijection between the mutation class of a quiver of type $D_n$
  and triangulations up to rotation and inverting tags} 

Here we adapt the methods and ideas of \cite{t} to obtain a
bijection between the mutation class of a quiver of type $D_n$ and
the set of triangulations of a punctured $n$-gon up to rotations and
inversion of tags. See also \cite{glz}.

\begin{figure}[htp]
  \begin{center}
    \includegraphics[width=4.0cm]{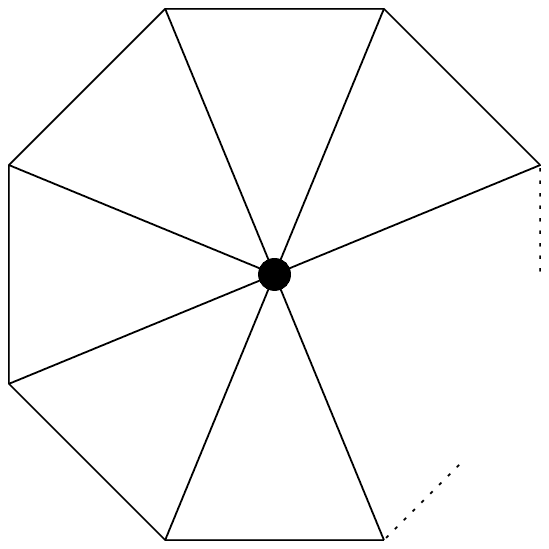}
  \end{center}\caption{Triangulation $R_n$ corresponding to the quiver
    consisting of the oriented cycle of length $n$.}
  \label{fignotoriented} 
  \end{figure}

We say that a diagonal from $a$
to $b$ is \textit{close to the border} if $|\delta(a,b)|=3$. For a
quiver $Q_\Delta$ corresponding to a triangulation $\Delta$, we will
always denote by $v_\alpha$ the vertex in $Q_\Delta$ corresponding
to the diagonal $\alpha$. From now on we let $n \geq 5$. Let us denote by
$S_n$ the triangulation of $\mathcal{P}_n$ shown in Figure
\ref{fignotoriented}. Note that this triangulation and the triangulation $S_n$ with all tags inverted are the only triangulations that correspond to the quiver consisting of the
oriented cycle of length $n$, $Q_n$.

\begin{lem}\label{exist diagonal}
Let $\Delta$ be a triangulation of $\mathcal{P}_n$, with $\Delta \neq
S_n$. Then there exists a diagonal in $\Delta$ which is close to the
border. 
\end{lem}
\begin{proof}
Let $\Delta$ be a triangulation of $\mathcal{P}_n$. If $\Delta$ is not
$S_n$, then there is at least one diagonal $\alpha$ which
connects two vertices on the border. See Figure
\ref{figdividedpolygon}.    

  \begin{figure}[h]
  \begin{center}
    \includegraphics[width=5cm]{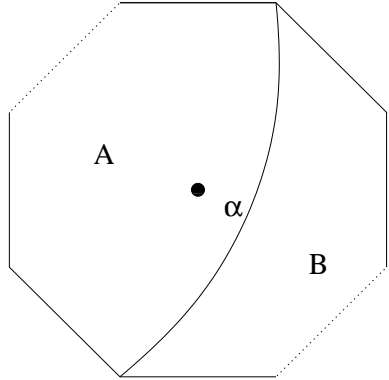}
  \end{center}
  \caption{The diagonal $\alpha$ divides the polygon into a punctured and a non-punctured surface.}
  \label{figdividedpolygon}
  \end{figure}

Consider the non-punctured surface $B$ determined by this diagonal. If
$\alpha$ is not close to the border, there exist a diagonal that
divides the surface $B$ into two smaller surfaces. By induction, there
exists a diagonal close to the border.  
\end{proof}

\begin{lem}\label{close to the border sink source cycle}
If a diagonal $\alpha$ of a triangulation $\Delta$ is close to the
border, then the corresponding vertex $v_{\alpha}$ in
$\epsilon_n(\Delta)=Q_{\Delta}$ is either a source, a sink or lies on
an oriented cycle of length 3. 
\end{lem}
\begin{proof}
Suppose $\alpha$ is a diagonal close to the
border. We have to consider the eight cases shown in Figure
\ref{figsinksourcecycleD}. In the first picture in Figure
\ref{figsinksourcecycleD}, $\alpha$ corresponds to a source since no other vertex
except $v_\beta$ can be adjacent to $v_\alpha$, or else the
corresponding diagonal would cross $\beta$. In the second picture
$\alpha$ corresponds to a sink. In picture three, four, five and six, there are
arrows between $v_\alpha$, $v_\beta$ and $v_{\beta'}$, and in the last
two pictures, there are arrows between $v_\alpha$, $v_\beta$ and
$v_\gamma$, so $v_\alpha$ lies on an oriented cycle of length $3$.
  \begin{figure}[htp]
  \begin{center}
    \includegraphics[width=3.0cm]{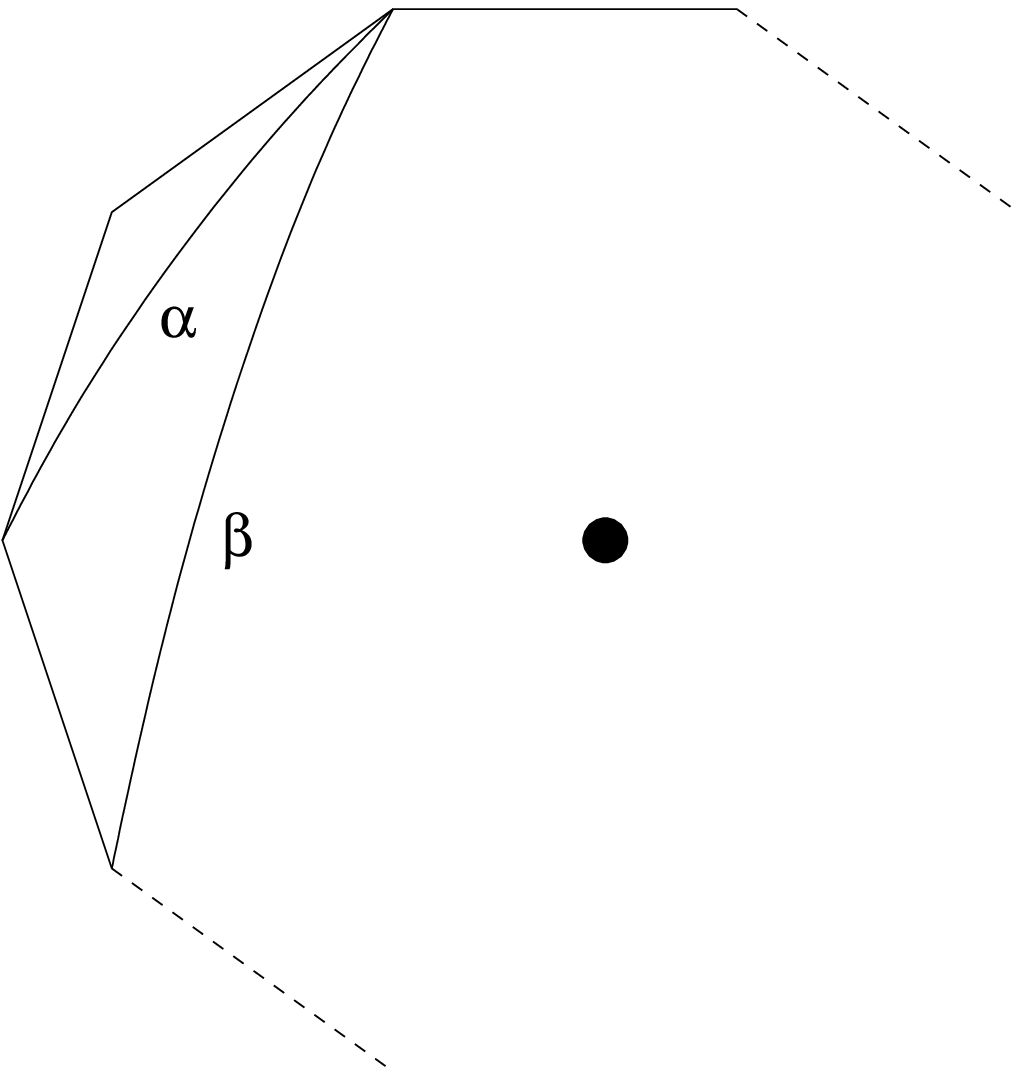}
    \includegraphics[width=3.0cm]{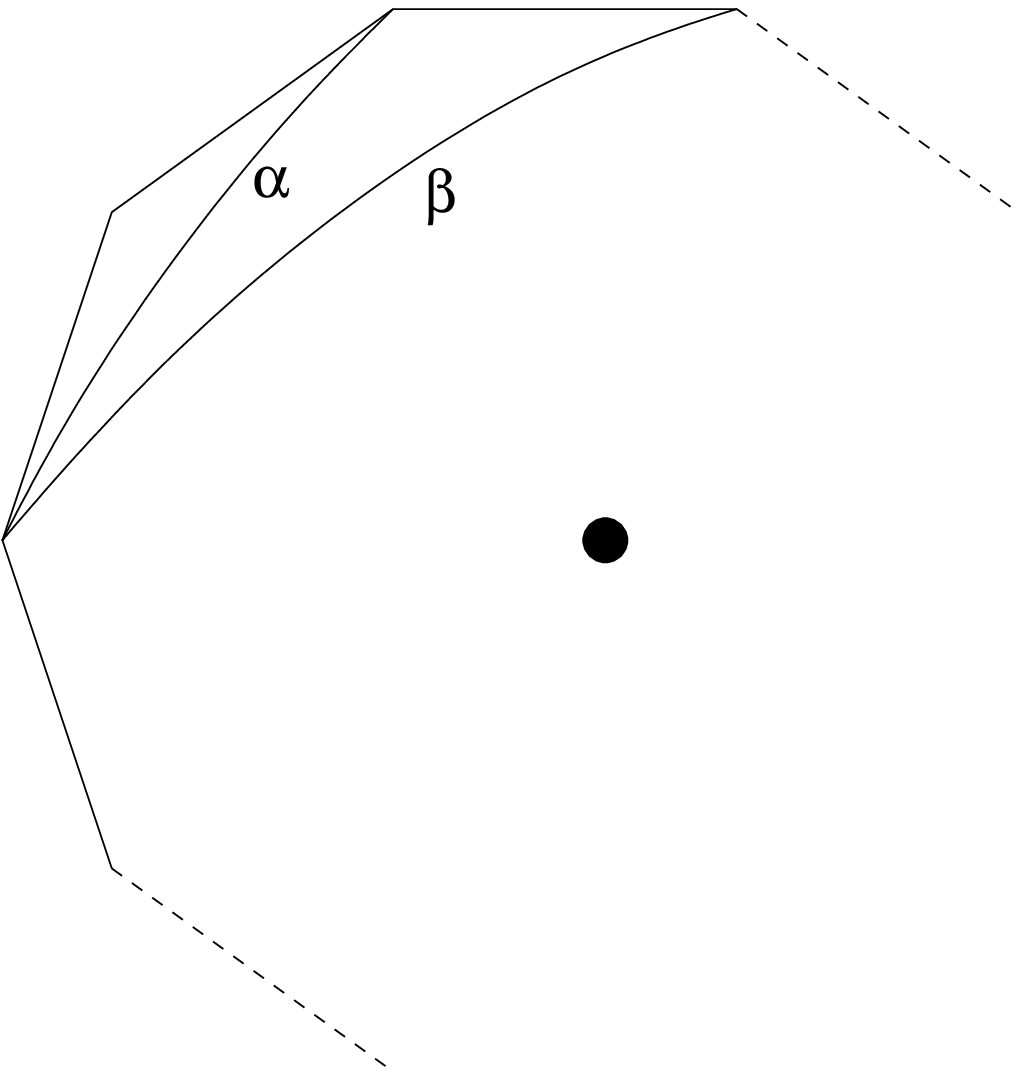}
    \includegraphics[width=3.0cm]{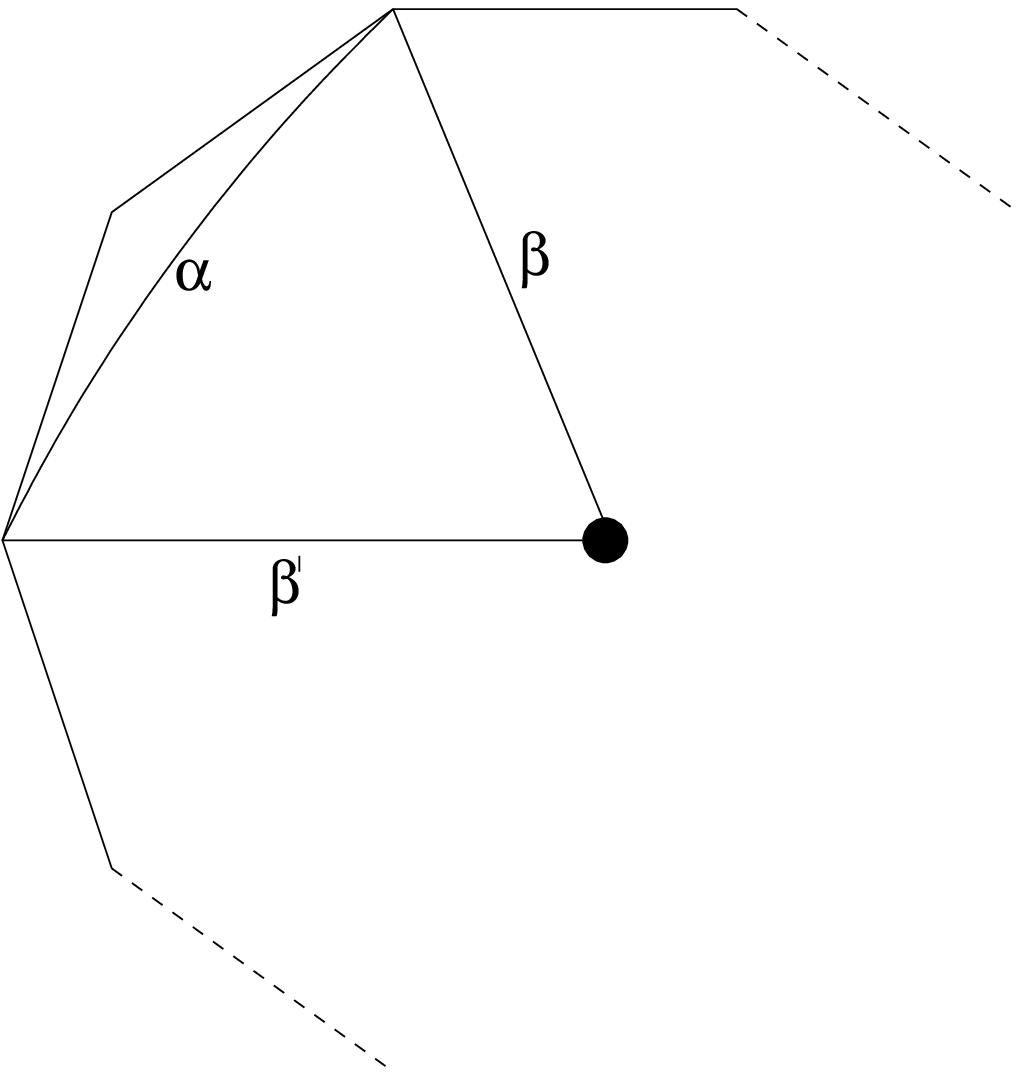}
    \includegraphics[width=3.0cm]{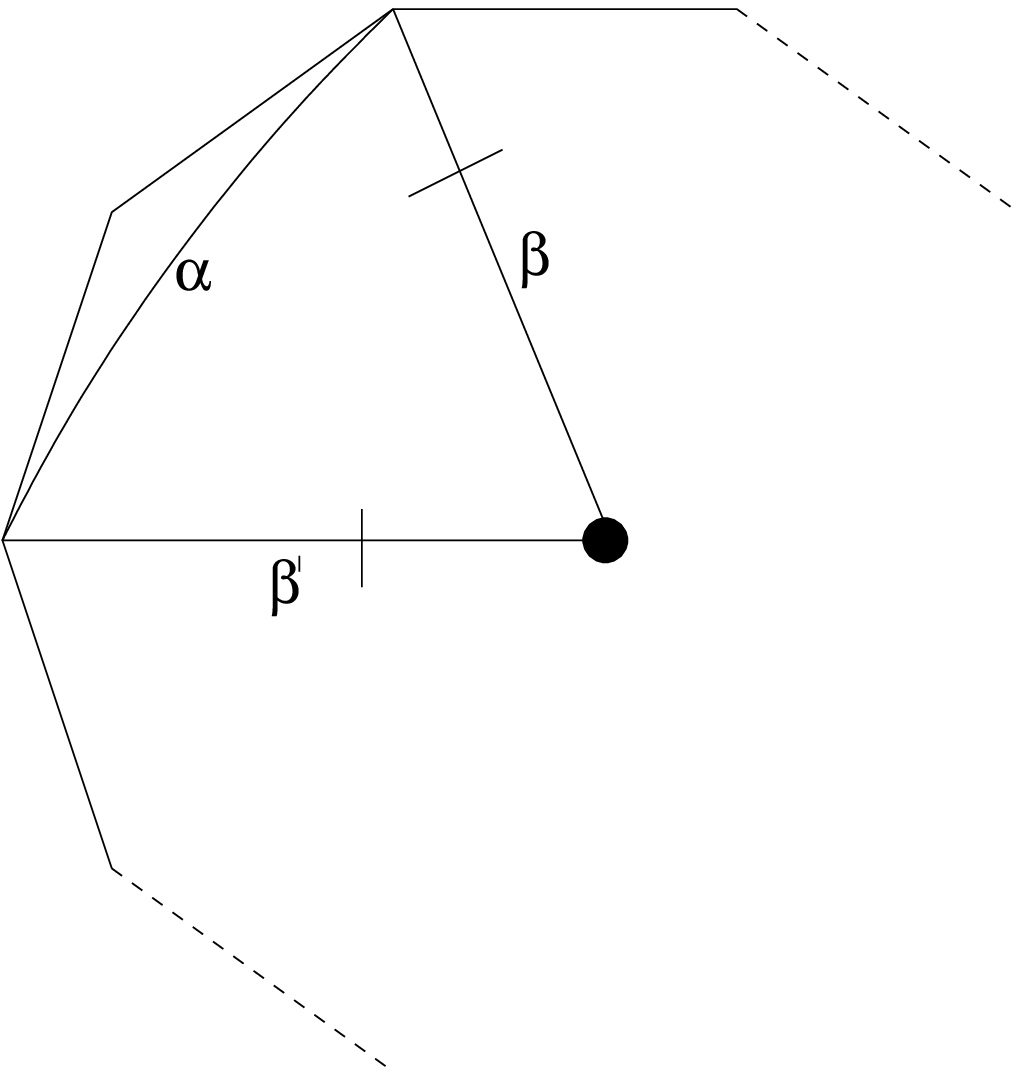}
    \includegraphics[width=3.0cm]{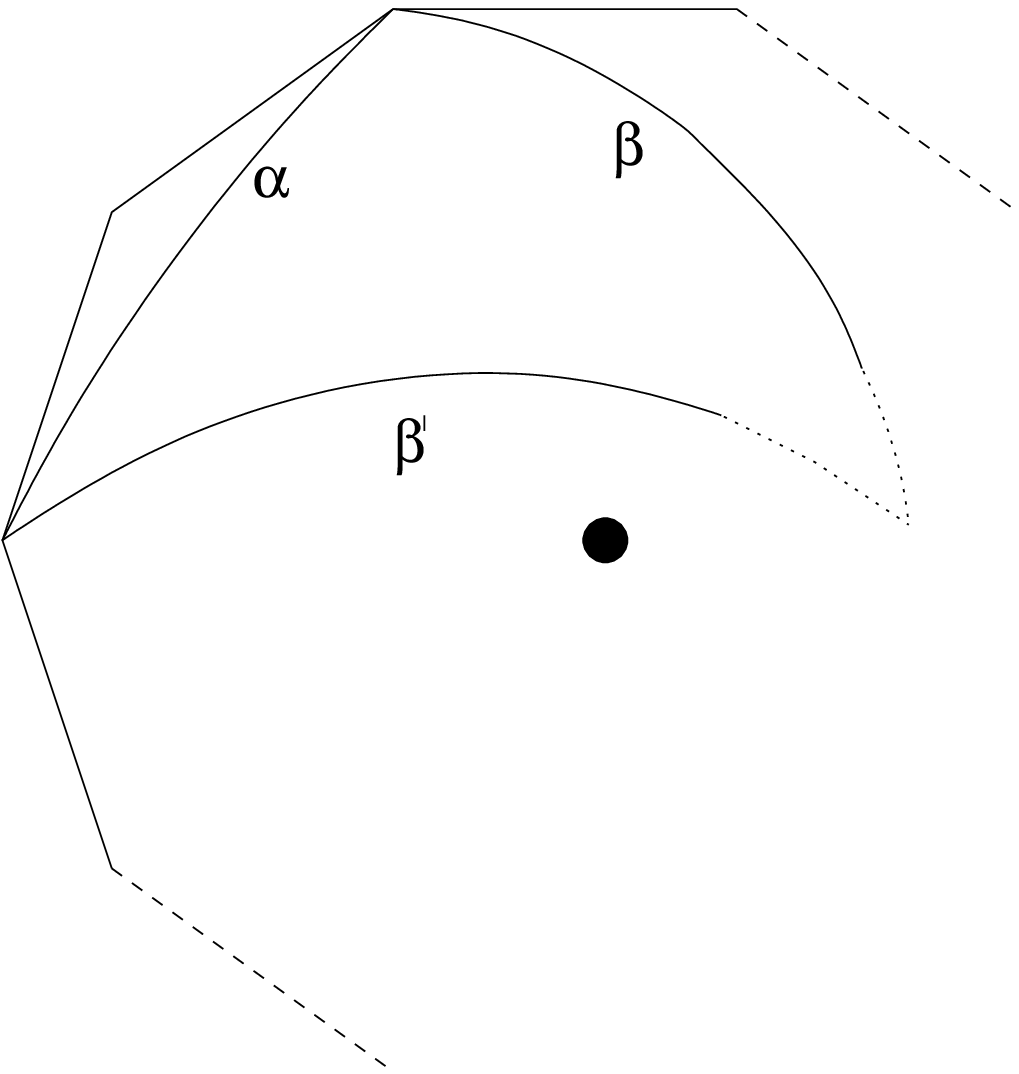}
    \includegraphics[width=3.0cm]{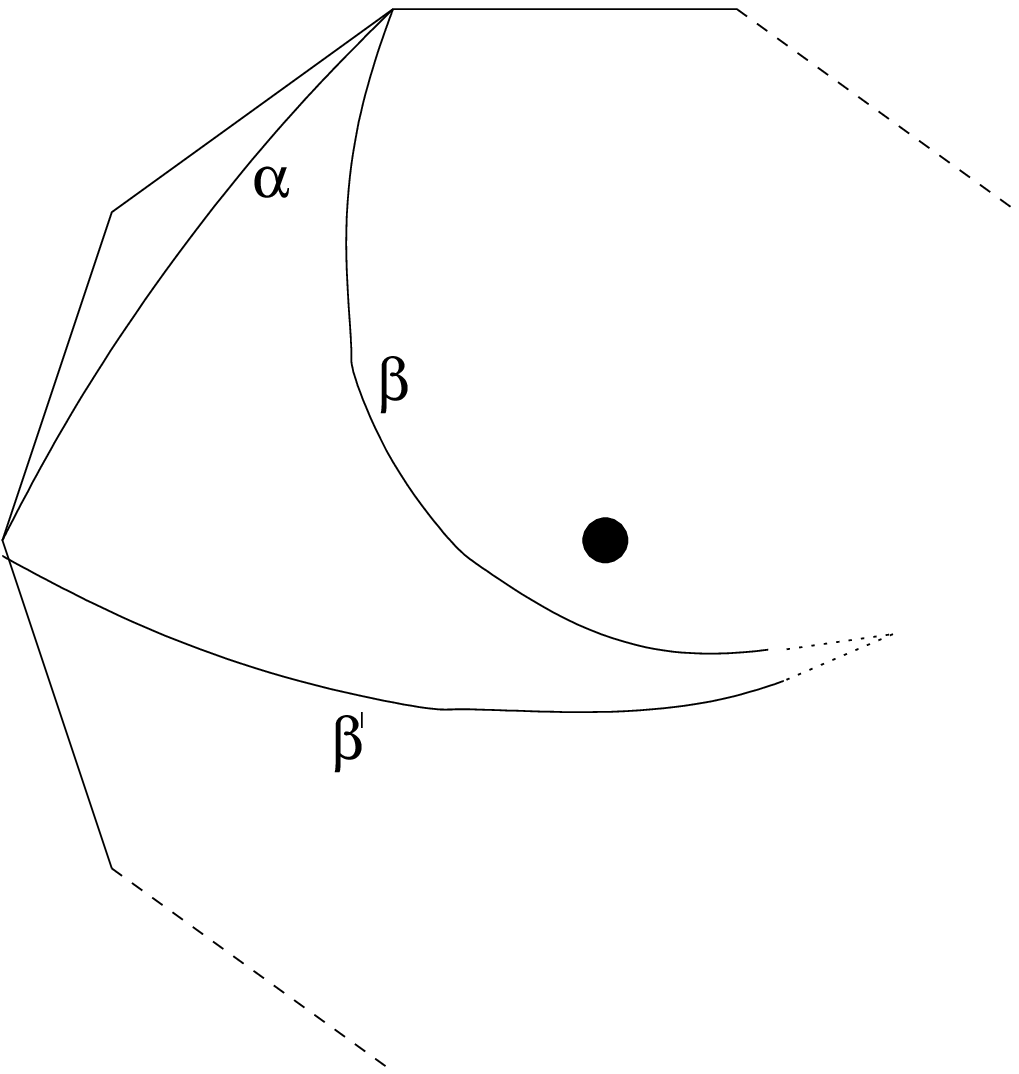}
    \includegraphics[width=3.0cm]{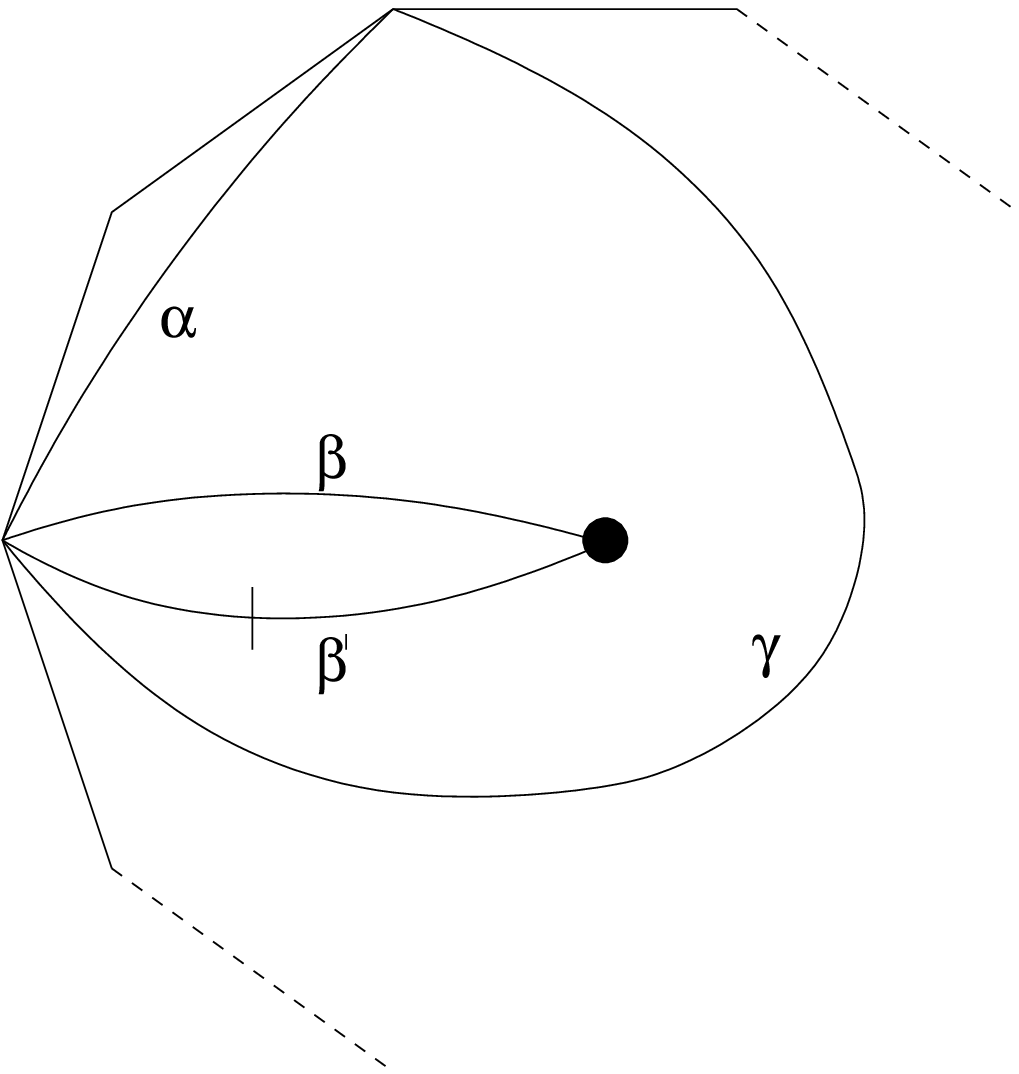}
    \includegraphics[width=3.0cm]{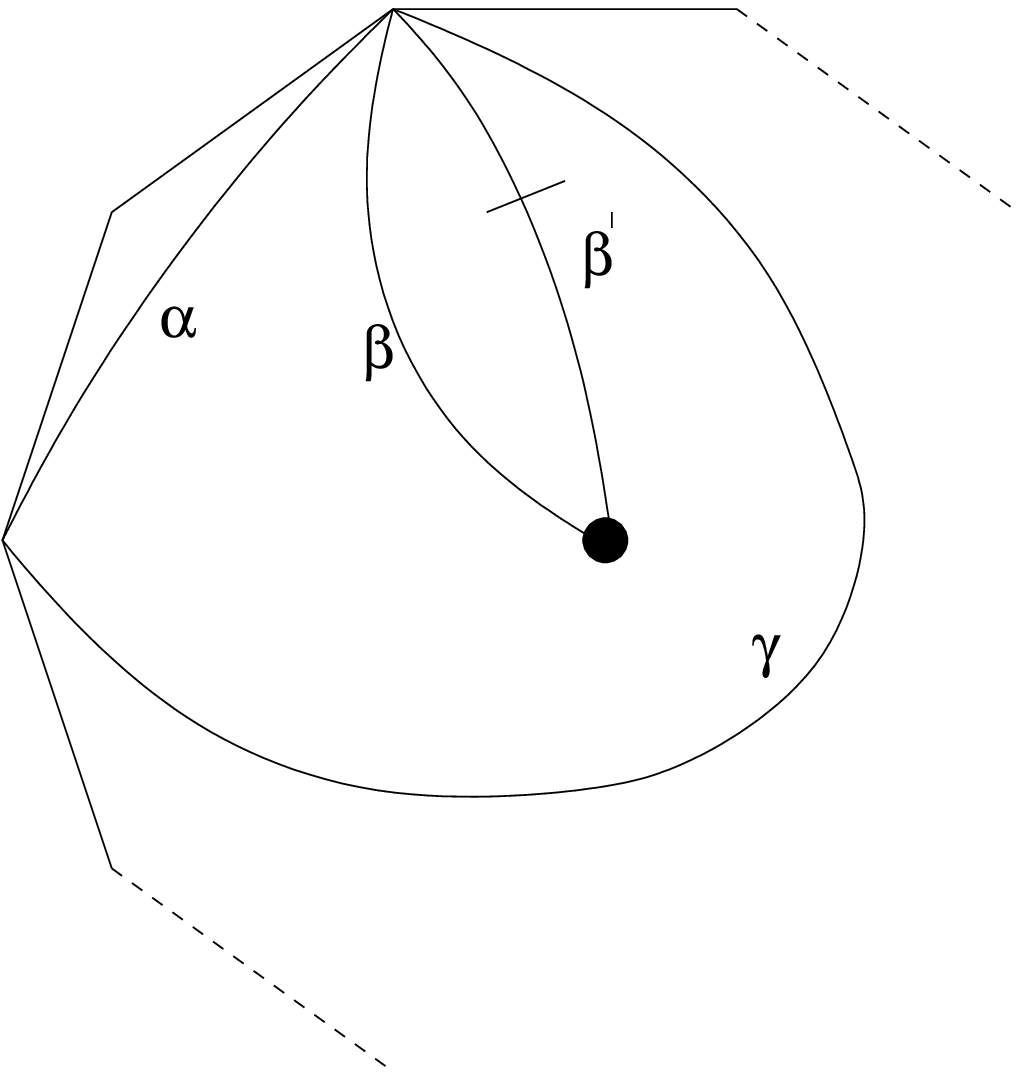}
  \end{center}\caption{See the proof of Lemma \ref{close to the border
      sink source cycle}}
  \label{figsinksourcecycleD} 
  \end{figure}
\end{proof}

\begin{figure}[h]
  \begin{center}
    \includegraphics[width=7cm]{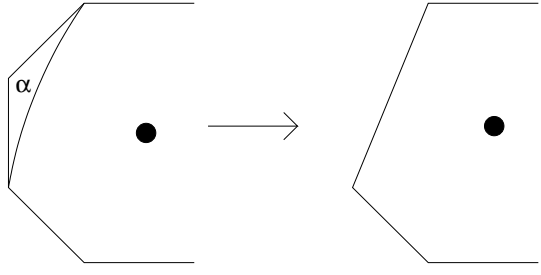}
  \end{center}
  \caption{Factoring out a diagonal close to the border.}
  \label{figfactoringD}
  \end{figure}

Let $\Delta$ be a triangulation of $\mathcal{P}_n$ and let
$\alpha$ be a diagonal close to the border. We define a triangulation
$\Delta / \alpha$ of $\mathcal{P}_n$ obtained from $\Delta$ by letting
$\alpha$ be a border edge and leaving all the other diagonals
unchanged. We write $\Delta/\alpha$ for the new triangulation
obtained and we say that we factor out $\alpha$. See Figure \ref{figfactoringD}. Note that this operation is
well-defined for each case in Figure \ref{figsinksourcecycleD}.    

\begin{lem}\label{is of type Dn}
Let $\Delta$ be a triangulation of $\mathcal{P}_n$, with $\Delta \neq S_n$ and let
$\epsilon_n(\Delta) = Q_{\Delta}$ be the corresponding quiver. If
$\alpha$ is a diagonal close to the border in $\Delta$, then the
quiver $Q_\Delta/v_\alpha$ obtained from $Q_\Delta$ by factoring out
the vertex $v_{\alpha}$ is connected and of type
$D_{n-1}$. Furthermore, we have that $\epsilon_{n-1}(\Delta / \alpha)
= Q_\Delta / v_\alpha$, when $\alpha$ is close to the border.  
\end{lem}
\begin{proof}
By Lemma \ref{close to the border sink source cycle} we have that $Q_\Delta / v_\alpha$
is connected. It is also straightforward to verify that
$\epsilon_{n-1}(\Delta/\alpha) = Q_\Delta /v_\alpha$ for each case, and hence
$Q_\Delta / v_\alpha$ is of type $D_{n-1}$ since $\Delta/\alpha$ is a
triangulation of $\mathcal{P}_{n-1}$.    
\end{proof}

Now we describe what happens when we factor out a vertex
corresponding to a diagonal not close to the border. We need to
consider two cases. We first deal with the case when $\alpha$ is a
diagonal not going between the puncture and the border.

\begin{lem}\label{not close to the border 1}
Let $\Delta$ be a triangulation and $\epsilon_n(\Delta)=Q_\Delta$. If we
factor out a vertex in $Q_{\Delta}$ corresponding to a diagonal that
is not close to the border and that is not a diagonal between the
puncture and the border, then the resulting quiver is disconnected.
\end{lem}
\begin{proof}
Let $\alpha$ be a diagonal not close to the border and not between the puncture and the border. Then the diagonal divides $\mathcal{P}_{n}$
into two surfaces $A$ and $B$. See Figure
\ref{figdividedpolygon}. Let $\beta$ be a diagonal in $A$ and $\beta'$
a diagonal in $B$. If $\beta$ and $\beta'$ would determine a common
triangle, the third diagonal would cross $\alpha$, hence there is no
arrow between the subquiver determined by $A$ and the subquiver
determined by $B$, except those passing through $v_\alpha$. It
follows that factoring out $v_\alpha$ disconnects the quiver.
\end{proof}

\begin{figure}[htp]
  \begin{center}
    \includegraphics[width=12.5cm]{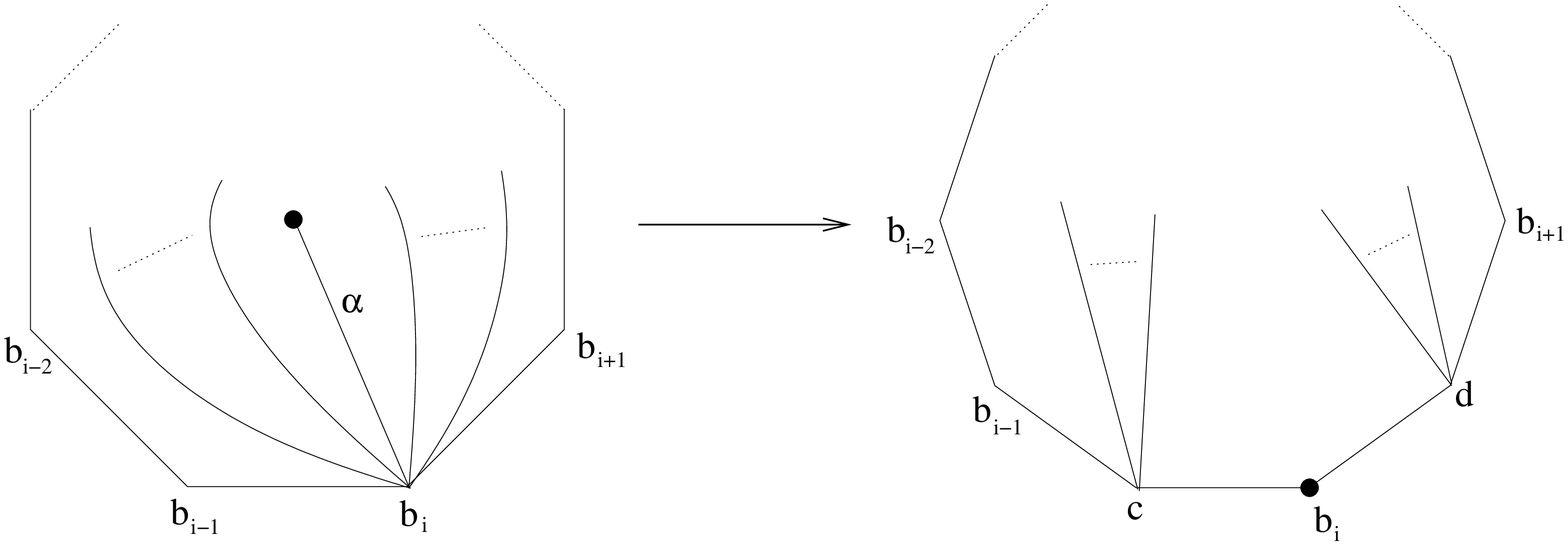}
    \includegraphics[width=12.5cm]{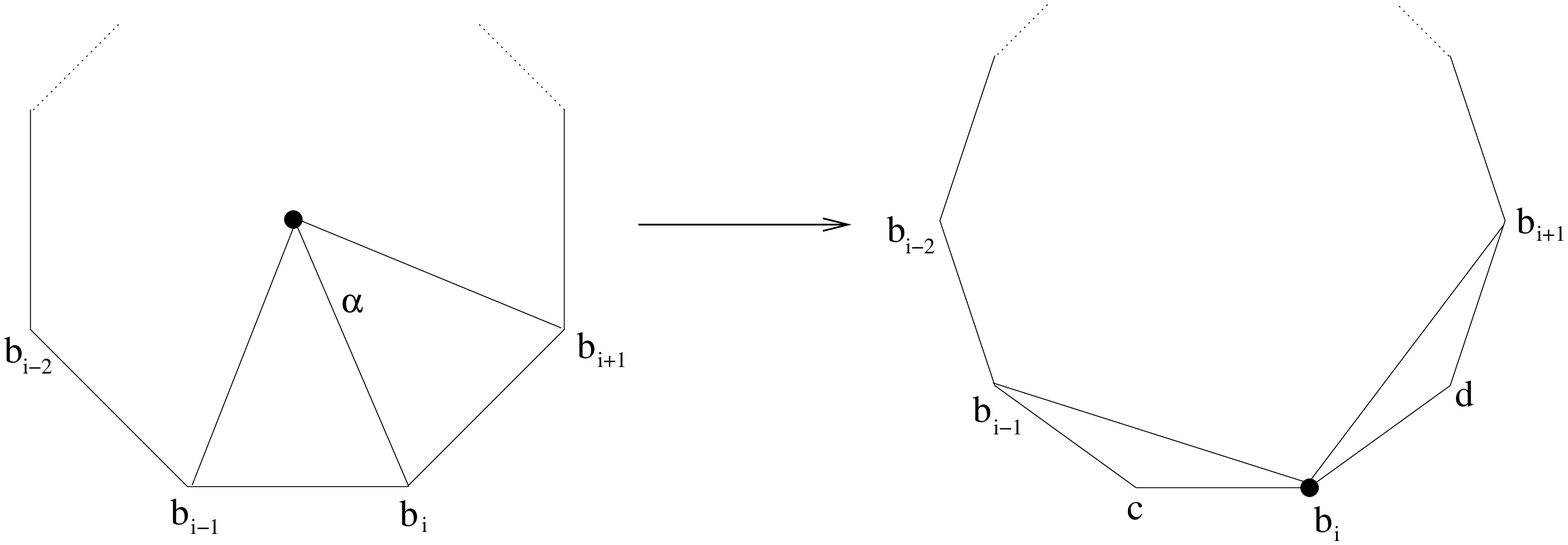}
  \end{center}\caption{Factoring out a diagonal from the puncture to
    the border.}
  \label{figfactpunctD}
  \end{figure}

Let $\Delta$ be a triangulation of $\mathcal{P}_n$ and let $\alpha$ be
a diagonal between the puncture and a vertex $b_i$ on the border of
the polygon. We want to understand the effect of factoring out
$v_\alpha$ (see Figure \ref{figfactpunctD}). In $\mathcal{P}_n$, create a new vertex $c$
between $b_{i-1}$ and $b_i$ and a new vertex $d$ between $b_i$ and
$b_{i+1}$, such that we obtain a $(n+2)$-polygon. Let all
diagonals that started in $b_i$ now start in $d$ and all diagonals
ending in $b_i$ now end in $c$. Remove the diagonal $\alpha$ and
identify the puncture with the vertex $b_i$. If there were two
diagonals between the puncture and $b_i$, remove both and draw a diagonal
from $c$ to $d$. Leave all the other diagonals unchanged. We will see
that this is a triangulation of the non-punctured $(n+2)$-polygon in the
next lemma.   

Recall that $\gamma_n$ is the function from the set of all
triangulations of the regular $(n+3)$-gon to the mutation
class of $A_n$, defined in Section 2. We have the following.

\begin{lem}\label{not close to the border 2}
Let $\Delta$ be a triangulation and $\epsilon_n(\Delta)=Q_\Delta$. If
$\alpha$ is a diagonal between the puncture and the border, then the
quiver $Q_\Delta / v_\alpha$ obtained from $Q_\Delta$ by factoring out
$v_\alpha$ is connected and of type $A_{n-1}$. Furthermore, we have
that $\gamma_{n+2}(\Delta/\alpha) = Q_\Delta / v_\alpha$ when $\alpha$
is a diagonal between the puncture and a vertex on the border. 
\end{lem}
\begin{proof}
It is clear that $\Delta / \alpha$ has $n-1$ diagonals and
that no diagonals cross. This means that the new triangulation is a
triangulation of the $(n+2)$ polygon without a puncture. We want to
show that all triangles are preserved by factoring out a diagonal as
described above and hence we will have that $\gamma_{n+2}(\Delta /
\alpha) = Q_\Delta / v_\alpha$, and that $Q_\Delta / v_\alpha$ is of
type $A_{n-1}$.    

First suppose that there is only one diagonal from the puncture to the
vertex $b_i$ (see Figure \ref{figfactpunctD}). Then it is easy to see
that all triangles are preserved. Next, suppose there are two
diagonals $\alpha$ and $\beta$ from the puncture to $b_i$. In this
case we add a new diagonal $\beta'$ between $b_{i-1}$ and $b_{i+1}$
and remove $\alpha$ and $\beta$. Then the diagonals bounding a common
triangle with $\beta$ before factoring out $\alpha$ will bound a
common triangle with $\beta'$ after factoring out $\alpha$. 
\end{proof}

Summarizing, we get the following Proposition.

\begin{prop}\label{if and only if}
Let $\Delta$ be a triangulation and let $\epsilon_n(\Delta)=Q_\Delta$ be
the corresponding quiver. Then $\epsilon_{n-1}(\Delta / \alpha) = Q_\Delta /
v_\alpha$ is of type $D_{n-1}$ if and only if the corresponding
diagonal $\alpha$ is close to the border. 
\end{prop}
\begin{proof}
From Lemma \ref{is of type Dn}, we have that if $\alpha$ is close to
the border, then $Q_\Delta / v_\alpha$ is of type $D_{n-1}$. If
$\alpha$ is not close to the border, we have by Lemma \ref{not close
  to the border 1} and Lemma \ref{not close to the border 2} that
$Q_\Delta / v_\alpha$ is either disconnected or of type $A_{n-1}$.  
\end{proof}

\begin{figure}[htp]
  \begin{center}
    \includegraphics[width=3.0cm]{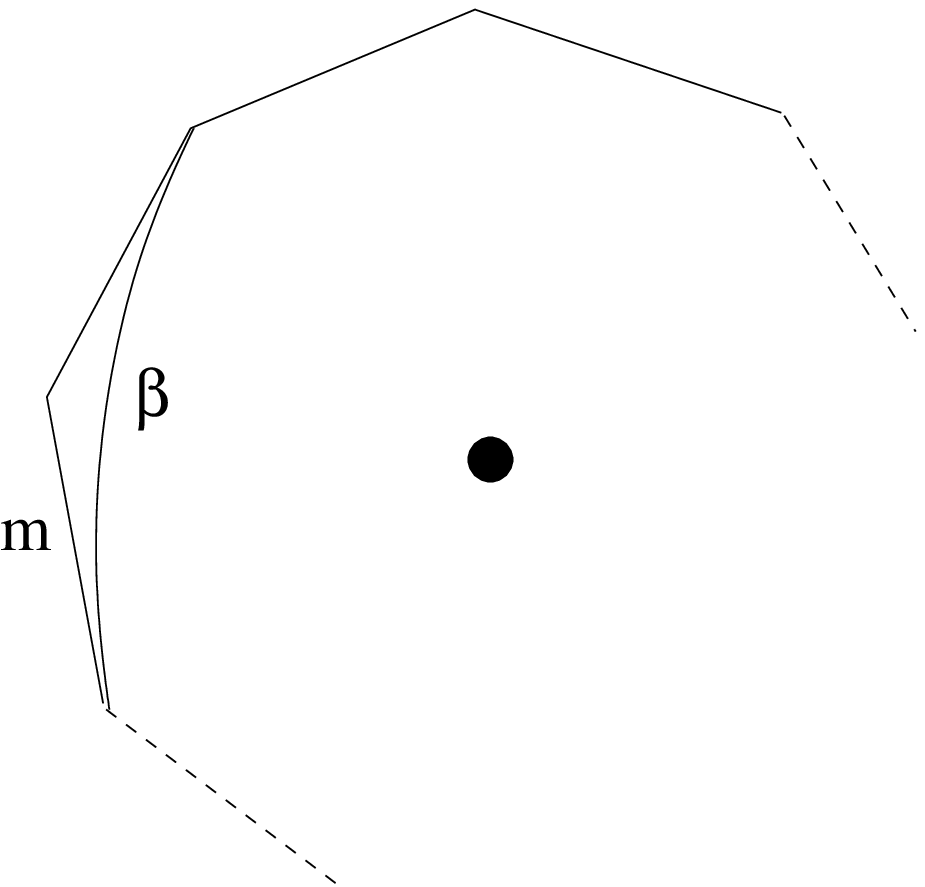}
    \includegraphics[width=3.0cm]{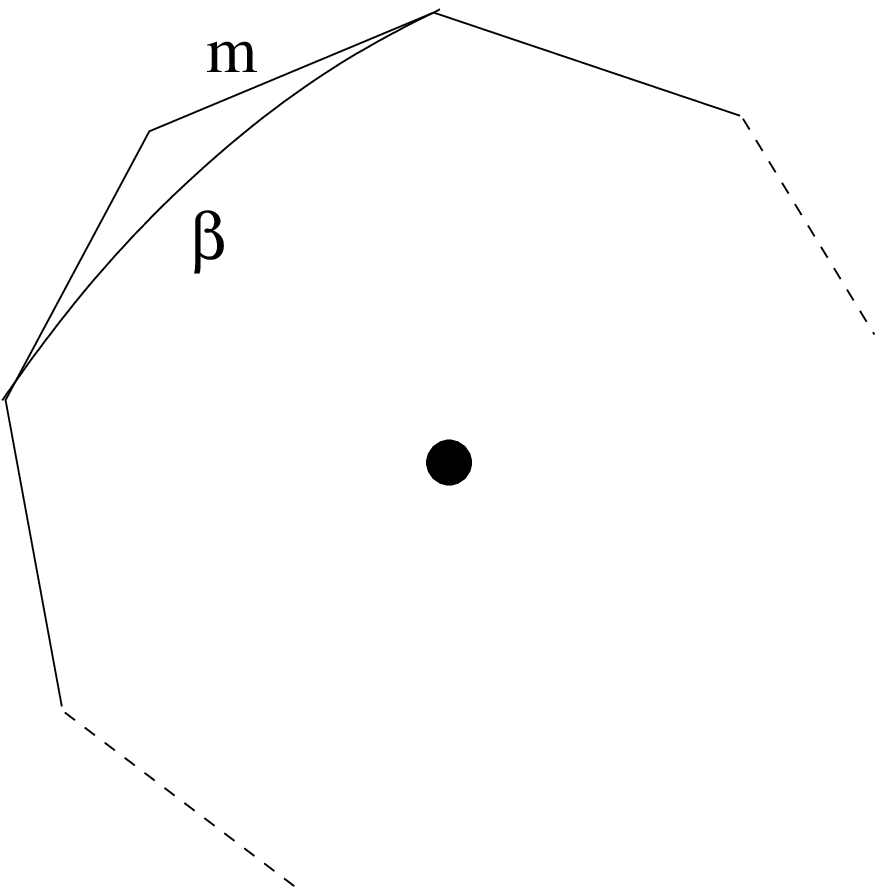}
    \includegraphics[width=3.0cm]{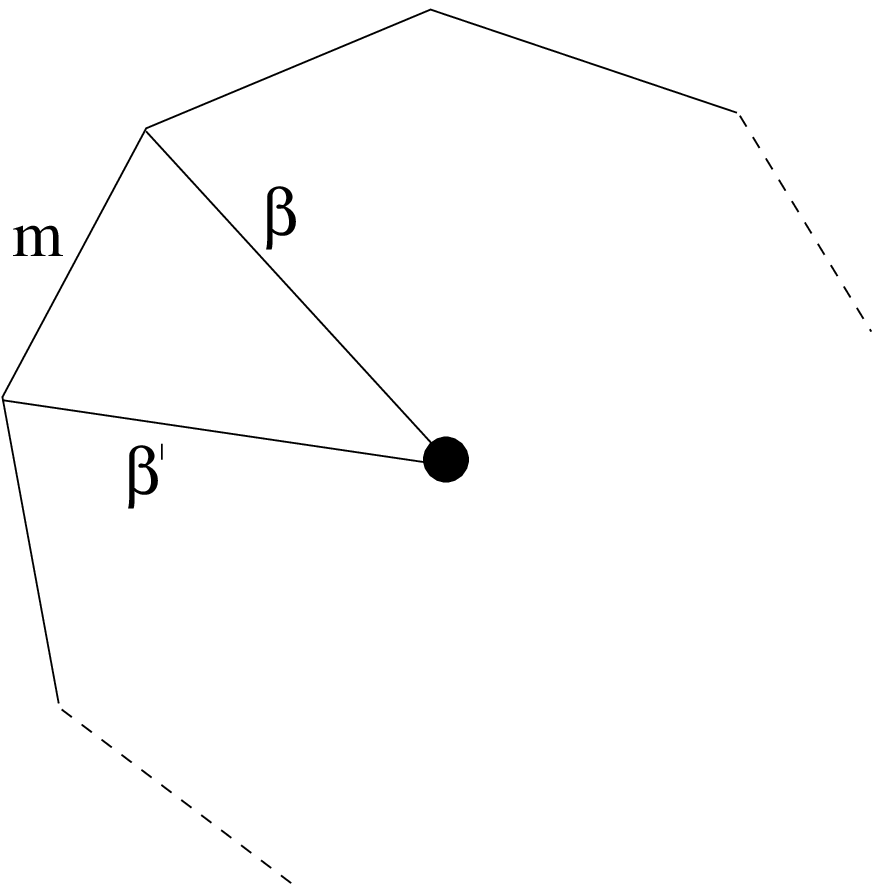}
    \includegraphics[width=3.0cm]{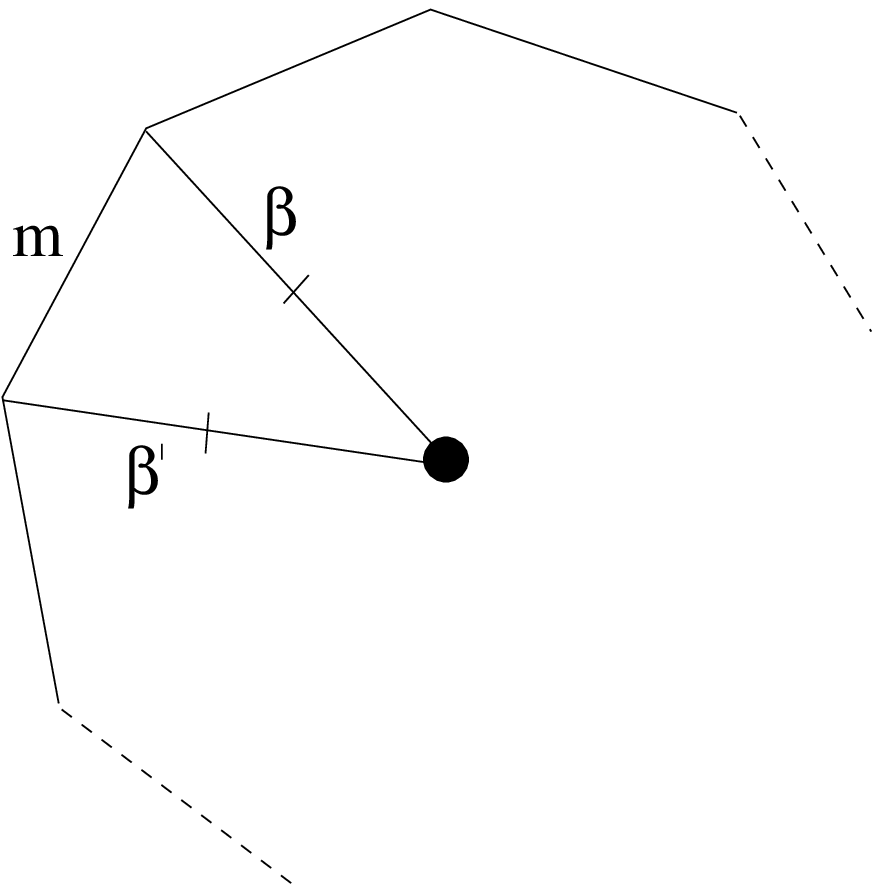}
    \includegraphics[width=3.0cm]{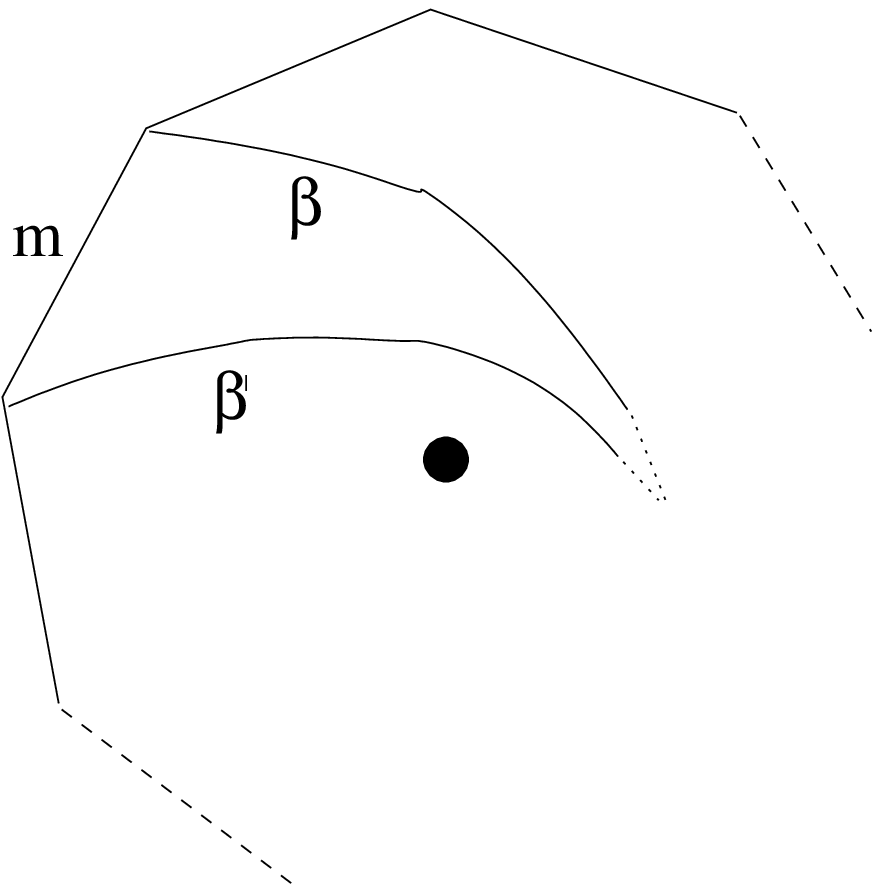}
    \includegraphics[width=3.0cm]{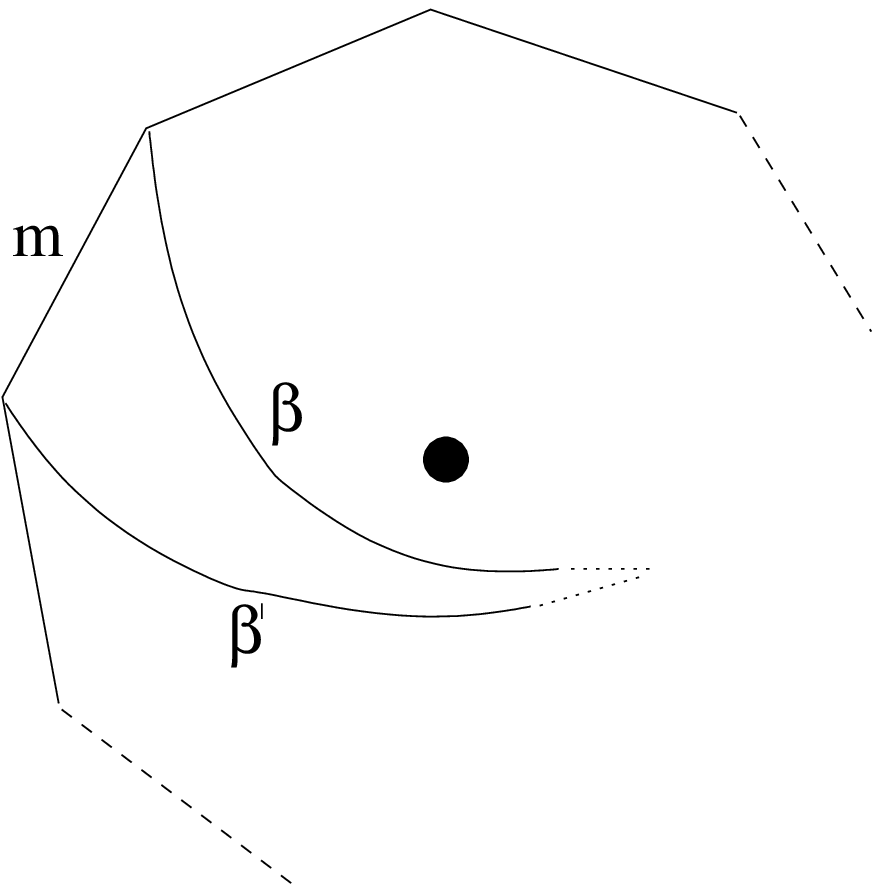}
    \includegraphics[width=3.0cm]{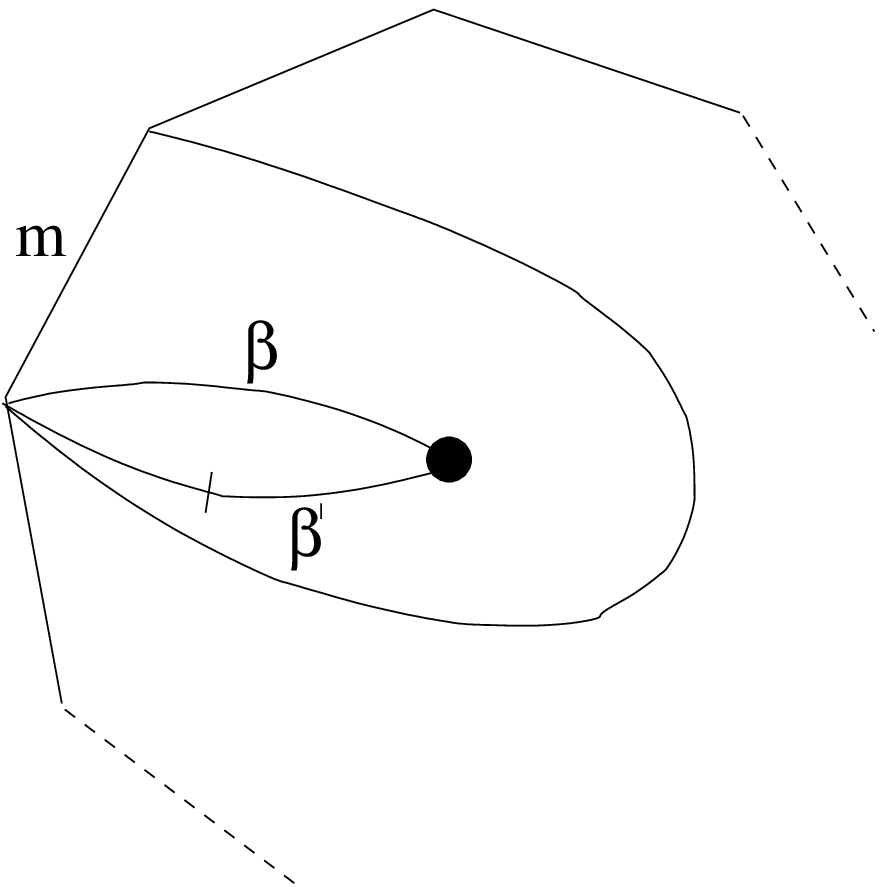}
    \includegraphics[width=3.0cm]{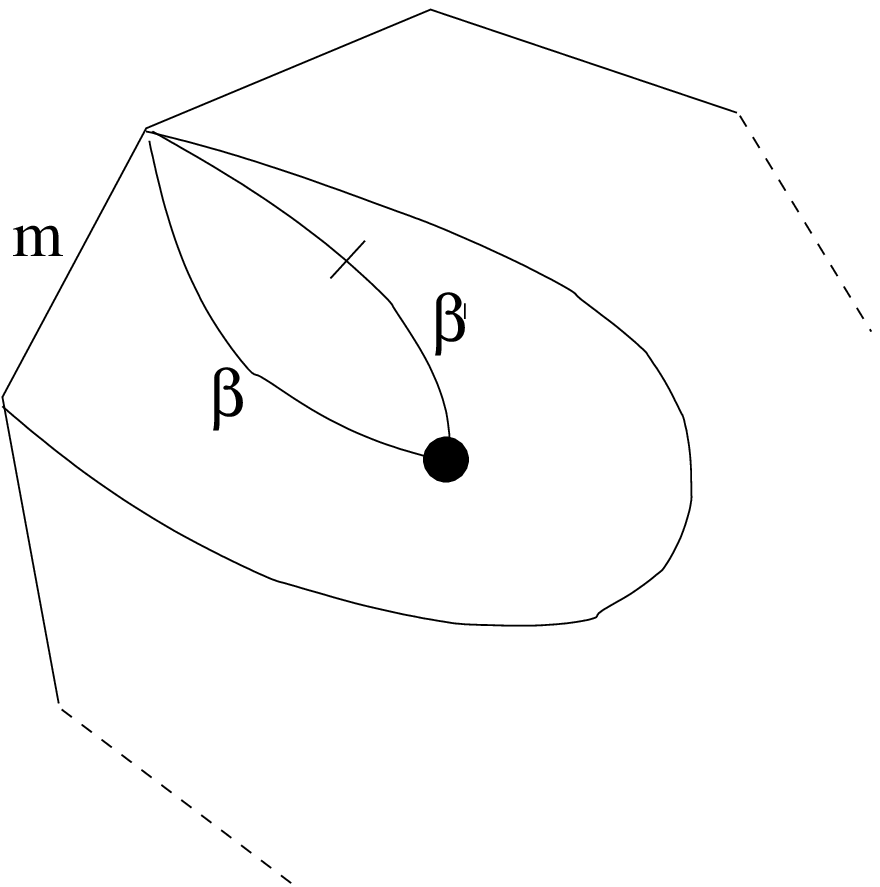}
  \end{center}\caption{Extension at $m$.}
  \label{figextendingD}
  \end{figure}

If $\Delta$ is a triangulation
of $\mathcal{P}_n$, we want to add a diagonal $\alpha$ and a vertex on
the polygon such that $\alpha$ is a diagonal close to the border and
such that $\Delta \cup \alpha$ is a triangulation of
$\mathcal{P}_{n+1}$. Consider any border edge $m$ on
$\mathcal{P}_n$. We consider the eight different cases for the
triangle containing $m$, as shown in Figure \ref{figextendingD}. We
can define the extension at $m$ for each case. See Figure
\ref{figsinksourcecycleD} for the corresponding extensions.   

For a given diagonal $\beta$, there are at most three
ways to extend the polygon with a diagonal $\alpha$ such that $\alpha$
is adjacent to $\beta$. These extensions
give non-isomorphic quivers, except when the triangulation is $S_n$. 

Combining Lemma \ref{exist diagonal} and Lemma \ref{is of type Dn}, we
get that for a quiver $Q$ which is not $Q_n$, there always exist a
vertex $v$ such that $Q'$ obtained from $Q$ by factoring out $v$ is
connected and a quiver of a cluster-tilted algebra of type
$D$. Furthermore, such a vertex must correspond to a diagonal close to
the border in any triangulation $\Delta$ such that $\epsilon_n(\Delta)
= Q_{\Delta}$.  

For a triangulation $\Delta$ of $\mathcal{P}_n$, let us denote by
$\Delta(i)$ the triangulation obtained from $\Delta$ by rotating $i$
steps in the clockwise direction. Also denote by $\Delta^{-1}$ the
triangulation obtained from $\Delta$ by inverting all tags. We define
an equivalence relation on $\mathcal{T}_n$ where we let $\Delta \sim 
\Delta(i)$ for all $i$ and $\Delta^{-1} \sim \Delta$. We define a new
function $\widetilde{\epsilon}_n:(\mathcal{T}_n / \! \! \sim) \rightarrow
\mathcal{M}_n$ induced from $\epsilon_n$. This is well-defined, and since $\epsilon_n$ is a surjection, we also have that
$\widetilde{\epsilon}_n$ is a surjection. We actually have the
following. 

\begin{thm}\label{maintheorem}
The function $\widetilde{\epsilon}_n:(\mathcal{T}_n / \! \! \sim) \rightarrow \mathcal{M}_n$
is bijective for all $n \geq 5$.
\end{thm}
\begin{proof}
We already know that $\widetilde{\epsilon}_n$ is surjective.

Suppose $\widetilde{\epsilon}_n(\Delta) =
\widetilde{\epsilon}_n(\Delta')$. We want to show that $\Delta =
\Delta'$ in $(\mathcal{T}_n / \! \! \sim)$ using induction.

It is straightforward to check that $\widetilde{\epsilon}_5:(\mathcal{T}_5 / \! \! \sim)
\rightarrow \mathcal{M}_5$ is injective. Suppose $\widetilde{\epsilon}_{n-1}:(\mathcal{T}_{n-1} / \! \! \sim)
\rightarrow \mathcal{M}_{n-1}$ is injective. Let $\alpha$ be a
diagonal close to the border in $\Delta$, with image $v_\alpha$ in
$Q$, where $Q$ is a representative for $\widetilde{\epsilon}_n(\Delta)$. Then
the diagonal $\alpha'$ in $\Delta'$ corresponding to $v_\alpha$ in $Q$
is also close to the border by Proposition \ref{if and only if}. We have
$\widetilde{\epsilon}_{n-1}(\Delta/\alpha)=\widetilde{\epsilon}_{n-1}(\Delta'/\alpha') =
Q/v_\alpha$, and hence by hypothesis, $\Delta/\alpha =
\Delta'/\alpha'$ in $(\mathcal{T}_n/ \! \! \sim)$.

We can obtain $\Delta$ and $\Delta'$ from $\Delta/\alpha =
\Delta'/\alpha'$ by extending the polygon at some border edge. Fix a
diagonal $\beta$ in $\Delta$ such that $v_\alpha$ and $v_\beta$ are
adjacent. This can be done since $Q$ is connected. Let $\beta'$ be the
diagonal in $\Delta'$ corresponding to $v_\beta$. By the above there
are at most three ways to extend $\Delta/\alpha$ such that the new
diagonal is adjacent to $\beta$. It is clear that these extensions
will be mapped by $\widetilde{\epsilon_n}$ to non-isomorphic quivers. Also
there are at most three ways to extend $\Delta'/\alpha'$ such that the
new diagonal is adjacent to $\beta'$, and all these extensions are
mapped to non-isomorphic quivers, thus $\Delta=\Delta'$ in
$(\mathcal{T}_n / \! \! \sim)$.
\end{proof}

\begin{cor} The number $d(n)$ of elements in the mutation class of any
  quiver of type $D_n$ is equal to the number of triangulations of the
  punctured regular $n$ polygon up to rotations and inverting all tags. 
\end{cor}

\section{Equivalences on the cluster category in the $D_n$ case}
Since the Auslander-Reiten translation $\tau$ is an equivalence, it is
clear that if $T$ is a cluster-tilting object in $\mathcal{C}$,
then the cluster-tilted algebras $\End_{\mathcal{C}}(T)^{\op}$ and
$\End_{\mathcal{C}}(\tau T)^{\op}$ are isomorphic. We know that $\tau$
corresponds to rotation of diagonals. In \cite{t} it was proven that
if $T$ and $T'$ are cluster-tilting objects in $\mathcal{C}$, then the
cluster-tilted algebras $\End_{\mathcal{C}}(T)^{\op}$ and
$\End_{\mathcal{C}}(T')^{\op}$ are isomorphic if and only if $T' = \tau^i T$
for an $i \in \mathbb{Z}$ in the $A_n$ case. 

Let $\alpha$ be a diagonal (indecomposable object in
$\mathcal{C}$). If $\alpha$ is a diagonal between the puncture and the
border, let $\alpha^{-1}$ denote the diagonal $\alpha$ with inverted
tag. We define

$$\mu \alpha = \begin{cases}\alpha^{-1}& \mbox{ if } \alpha \mbox{ is a
    diagonal between the puncture and the border,} \\ \alpha&
  \mbox{ otherwise.}\end{cases}$$

\begin{figure}[htp]
  \begin{center}
    \includegraphics[width=12.0cm]{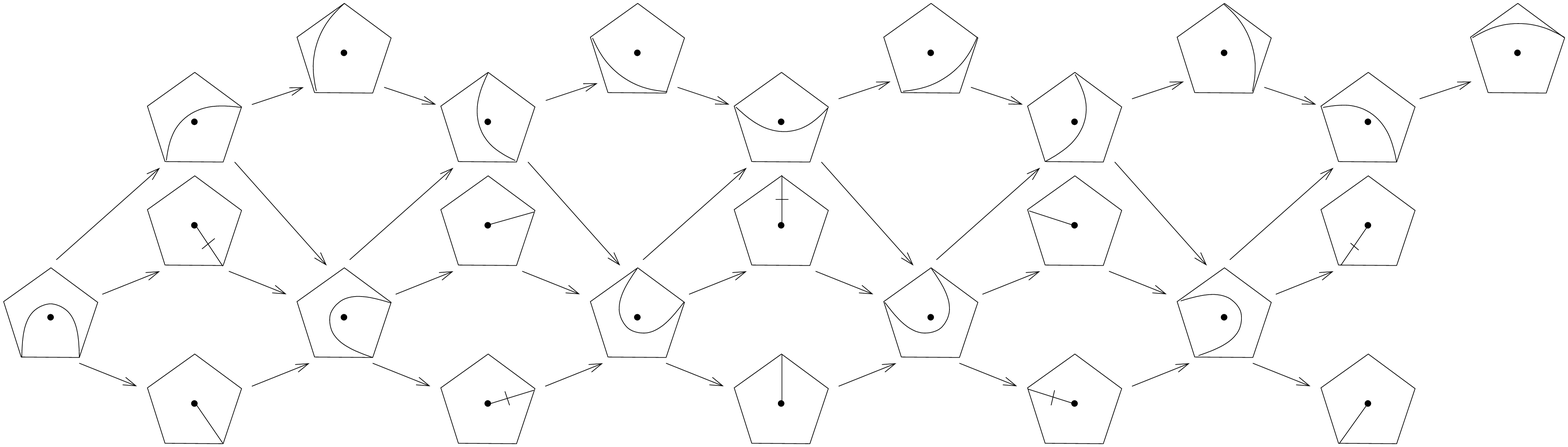}
  \end{center}\caption{AR quiver for the cluster category in the $D_5$ case.}
  \label{figARquiverD_5}
  \end{figure}

If $\alpha$ is not a diagonal between the puncture and the
border, then clearly $\tau^n \alpha = \alpha$. Now, let $\alpha$ be
a diagonal between the puncture and the border. Suppose $n$ is even. Then
it is clear from combinatorial reasons that $\tau^n \alpha = \alpha$
and that $\tau^i \alpha \neq \alpha^{-1}$ for any $i$. If $n$ is odd,
then $\tau^n \alpha = \alpha^{-1}$ and hence $\tau^{n} = \mu$. See
Figure \ref{figARquiverD_5} for an example of an AR-quiver in the
$D_5$ case.

\begin{thm}\label{maintheorem2}
Let $T$ and $T'$ be cluster-tilting objects in $\mathcal{C}$. Then the
cluster-tilted algebras $\End_{\mathcal{C}}(T)^{\op}$ and
$\End_{\mathcal{C}}(T')^{\op}$ are isomorphic if and only if $T' =
\mu^i \tau^j T$ $i,j \in \mathbb{Z}$.
\end{thm}
\begin{proof}
Let $\Delta$ be a triangulation corresponding to $T$ and $\Delta'$ a
triangulation corresponding to $T'$. If $T' \not\simeq \tau^i T$ for any
$i$, then $\Delta'$ is not obtained from $\Delta$ by a rotation. If
$T' \not\simeq \mu T$, then $\Delta \neq \Delta^{-1}$. It then follows
from Theorem \ref{maintheorem} that $\End_{\mathcal{C}}(T)^{\op}$ is not
isomorphic to $\End_{\mathcal{C}}(T')^{\op}$.
\end{proof}

It is clear that $\mu$ is an equivalence on the cluster category, since $\mu^2 =$ id.

\section{The number of triangulations of punctured polygons}
In this section we want to find an explicit formula for the number of
triangulations of punctured polygons up to rotation and tags. Let
$\mathcal{B}_n$ be the set of equivalence classes of trees such that 
\begin{itemize}
\item any full subtree not including the root is binary and every
  inner node has either two or no children, 
\item there are exactly $n$ leaves and
\item two trees are equivalent if one can be obtained from the other
  by rotating at the root. 
\end{itemize}

\begin{figure}[htp]
  \begin{center}
    \includegraphics[width=10cm]{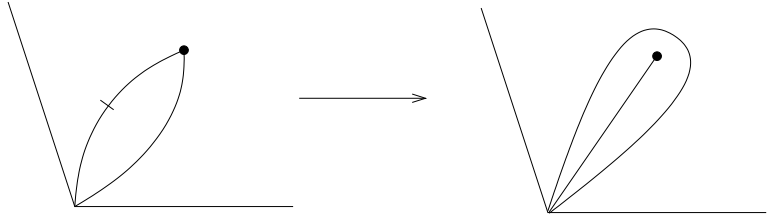}
  \end{center}\caption{Drawing tagged edges as loops}
  \label{figtagtoloop}
  \end{figure}

As before, let $\mathcal{T}_5 / \! \! \sim$ be the set of triangulations of the
punctured $n$-gon, where rotations and inverting tags gives equivalent
triangulations. In this section we will draw certain tagged edges as
loops. If there are two diagonals between the puncture and the same
vertex, we will draw one diagonal as a loop. See Figure
\ref{figtagtoloop}. 

\begin{figure}[htp]
  \begin{center}
    \includegraphics[width=12.5cm]{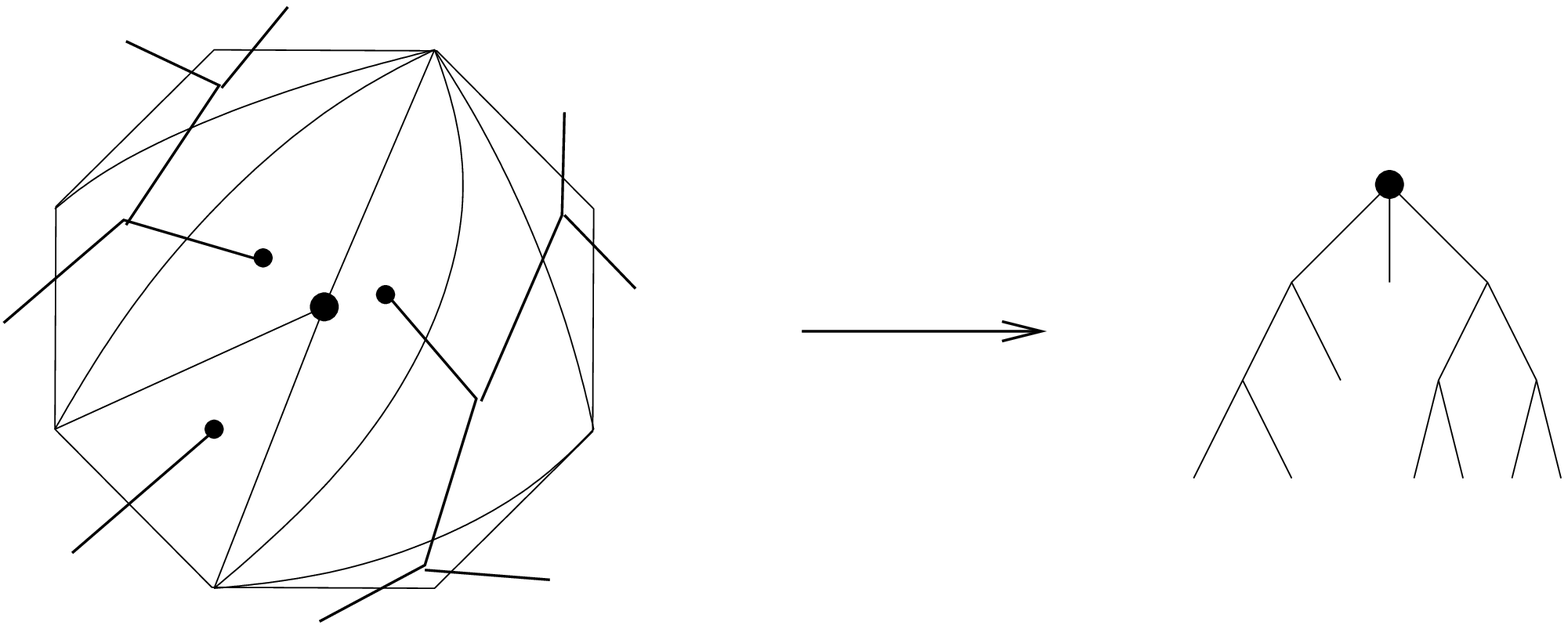}
    \includegraphics[width=12.5cm]{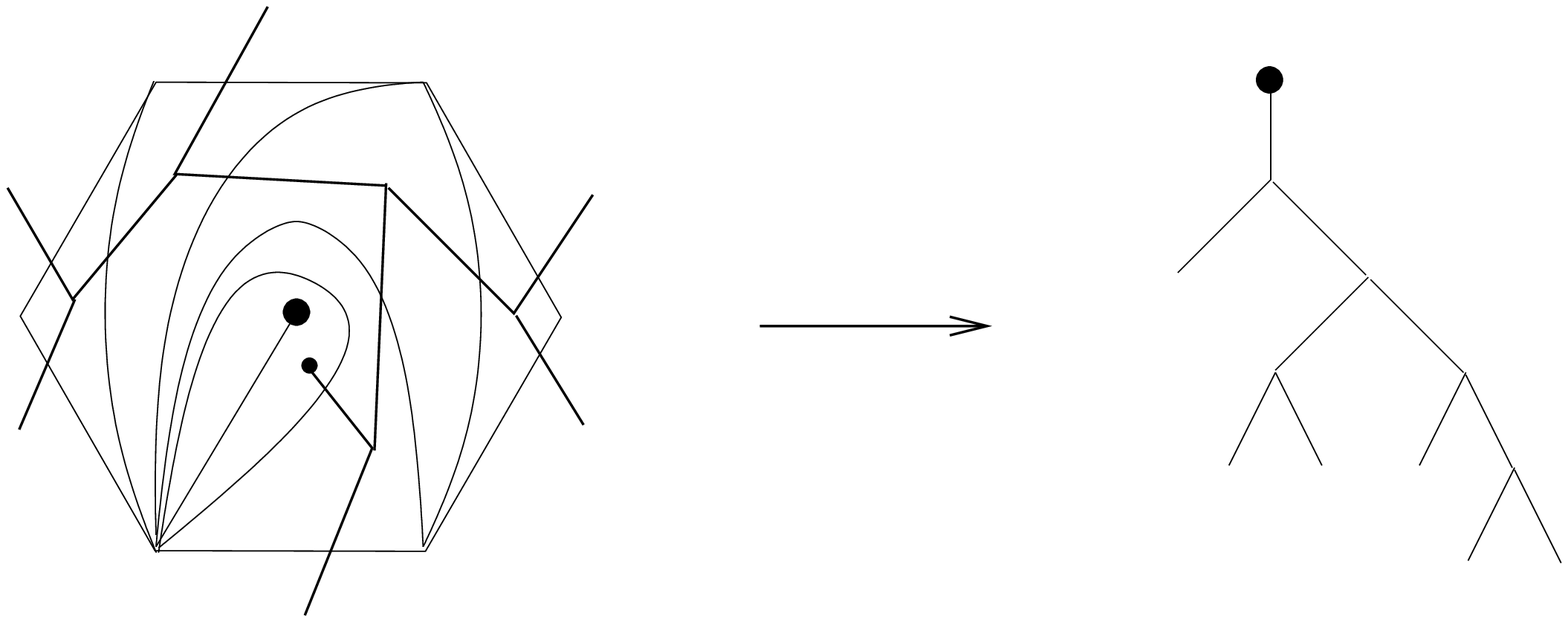}
  \end{center}\caption{Triangulation and corresponding tree.}
  \label{figtreitree}
  \end{figure}

We define a function $\sigma : \mathcal{T}_n / \! \! \sim \rightarrow
\mathcal{B}_n$ by assigning to a triangulation a tree. Let $\Delta$ be
a triangulation. We let $\sigma(\Delta)$ be the tree obtained in the
following way. Draw an edge between two triangles $E$ and $E'$ if they
are adjacent and their common diagonal is not a diagonal between the
puncture and the border. Note that a loop in this case is not an edge
between the puncture and the border. When a triangle $E$ contains one
or two border edges, also draw one or two edges from the vertex to the
outside of the polygon, crossing the border edges. These will be the
leaf edges. Then identify the vertices adjacent to the
puncture to be the root in the tree. See Figure \ref{figtreitree} for
some examples.   

It is clear that $\sigma$ is a well-defined function. Our aim is to
show that $\sigma$ is a bijection. 

\begin{figure}[htp]
  \begin{center}
    \includegraphics[width=3cm]{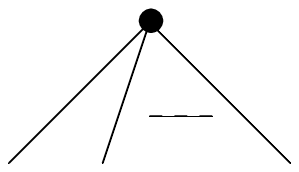}
  \end{center}\caption{The tree $R_n$ consisting of exactly $n$ edges from
the root.}
  \label{root tree}
  \end{figure}

Let the tree $R_n$ be the tree consisting of exactly $n$ edges from
the root, as shown in Figure \ref{root tree}. Note that this is the unique tree
which is the image of the triangulation $S_n$.

Now we want to define a function $\lambda:\mathcal{B}_n \rightarrow
\mathcal{T}_n / \! \! \sim$ and we will see that this is the inverse of $\sigma$. 

Given a tree $T$ with $n$ leaves, we will here describe
$\lambda(T)$. We know that an inner edge of a tree (an edge not
going to a leaf) corresponds to a diagonal $\alpha$ not going between
the puncture and the border.  

\begin{figure}[htp]
  \begin{center}
    \includegraphics[width=10.0cm]{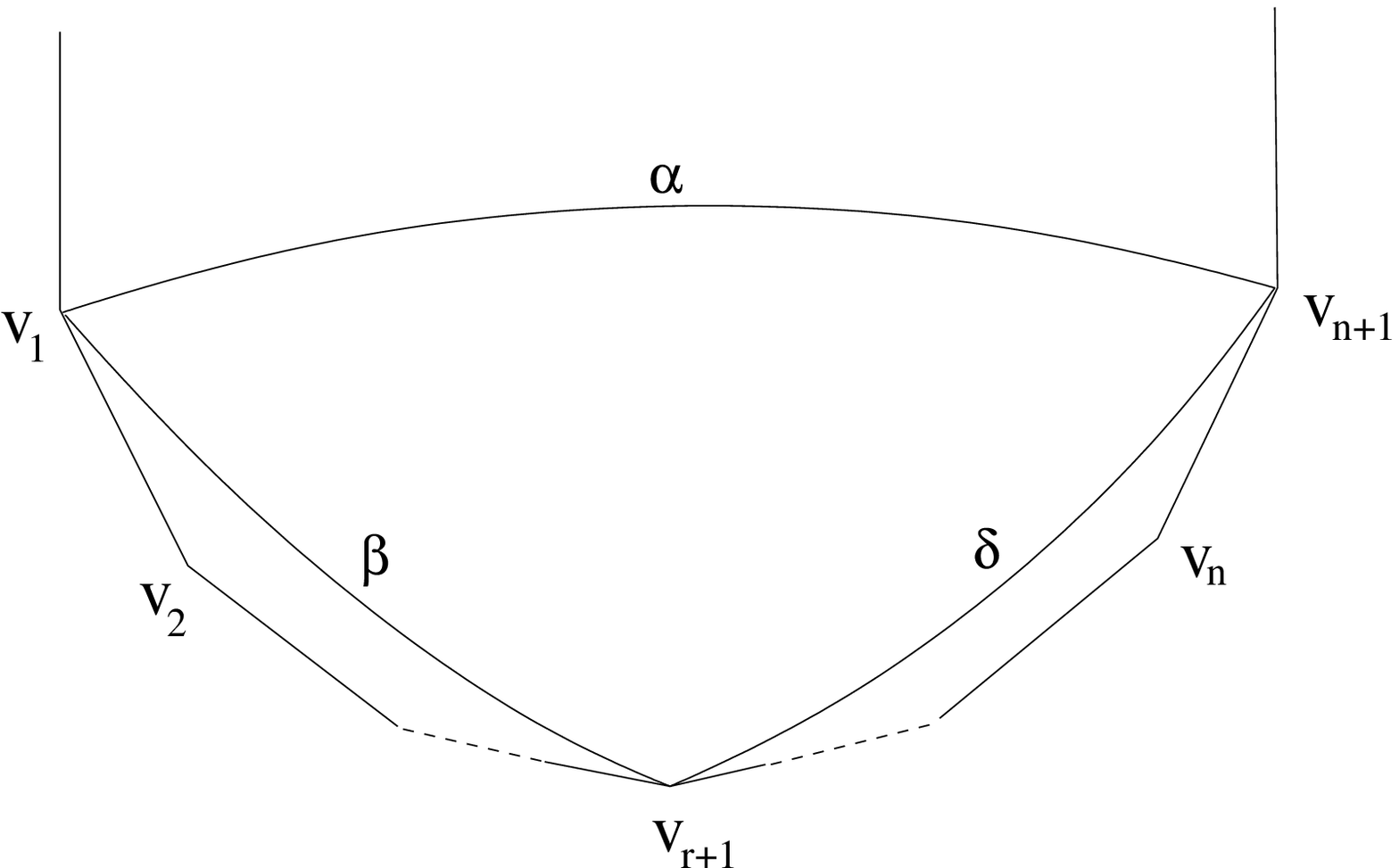}
  \end{center}\caption{}
  \label{correspondingpart2}
  \end{figure}

Suppose $\alpha$ is an inner
edge of $T$. Let $T'$ be the full subtree of $T$ with root ending in
$\alpha$. If $T'$ has $n$ leaves, we draw a segment of a polygon
consisting of $n$ border edges. See Figure
\ref{correspondingpart2}. Suppose the subtree to the left of the root
in $T'$ has $r \geq 2$ leaves. Then we draw a diagonal $\beta$ from
$v_1$ to $v_{r+1}$. If $r+1 \neq n$ we draw a diagonal $\delta$ from
$v_{r+1}$ to $v_{n+1}$. We can continue like this with $\beta$ and
$\delta$ until we made a complete triangulation of the segment of the
polygon, by induction.

Now, suppose $T$ has $k$ edges from the root, namely $t_1 , t_2 ,...,
t_k$. Suppose the full subtree with root ending in $t_i$ has $d_i$
leaves. Then $\sum_i d_i = n$. Draw a punctured polygon with $n$
border edges and draw $k$ diagonals between the puncture and vertices
on the border such that each segment has $d_i$ border edges in
anticlockwise direction. 

\begin{figure}[htp]
  \begin{center}
    \includegraphics[width=10.0cm]{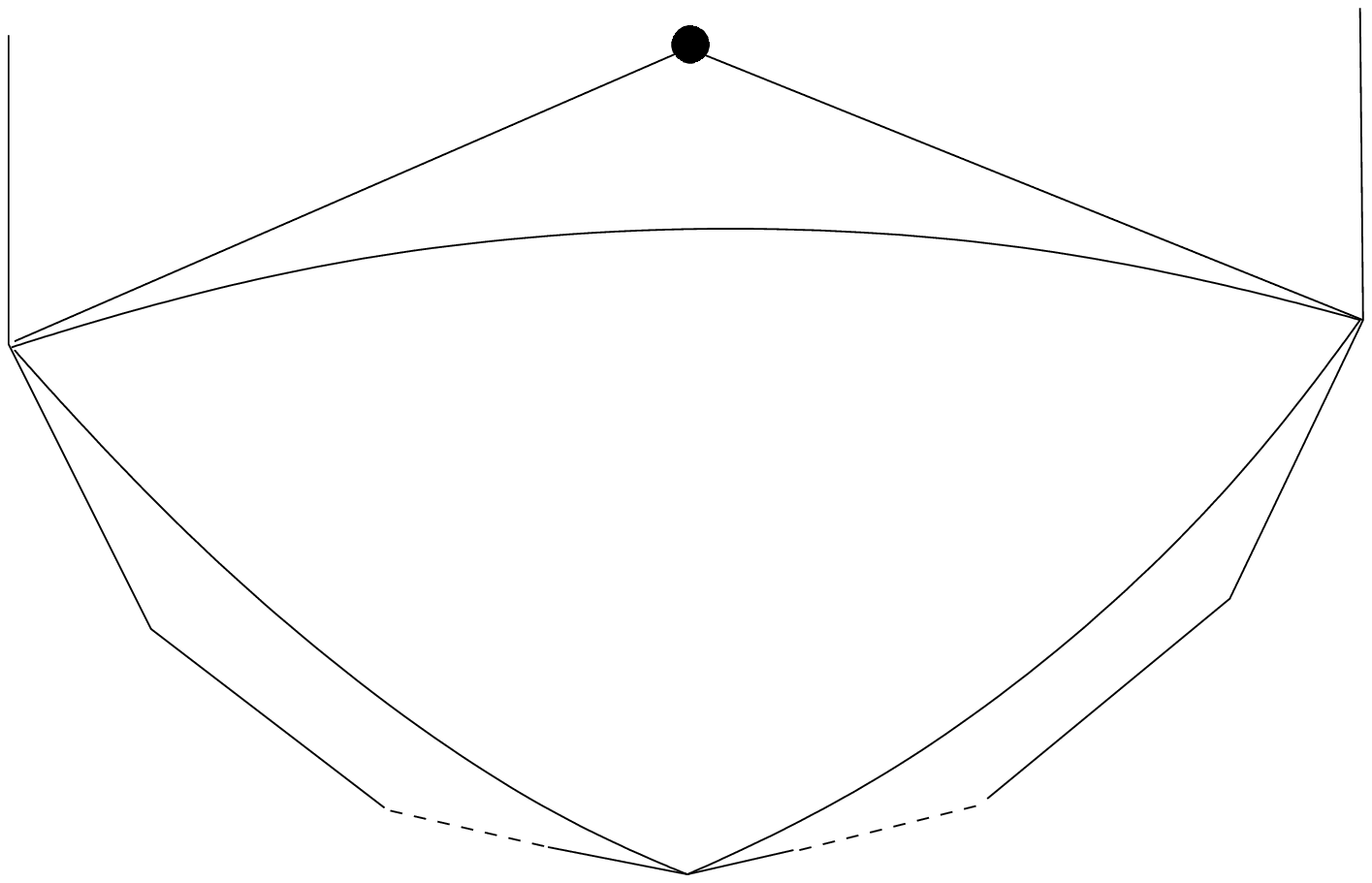}
  \end{center}\caption{}
  \label{correspondingpart}
  \end{figure}

For each segment defined by $t_i$, apply the procedure described above to
obtain a triangulation of the segment. See Figure
\ref{correspondingpart}.  

It is clear from the construction that $\lambda$ is the inverse of
$\sigma$, so we have the following.

\begin{thm}
$\sigma: \mathcal{T}_n \rightarrow \mathcal{B}_n$ is a bijection.
\end{thm}

The number of rooted planar trees with $n+1$ nodes where rotating at
the root gives equivalent trees, is given by the formula

\[\sum_{d|n} \phi(n/d)\binom{2d}{d}/(2n)\]
where $\phi$ is the Euler function (see \cite{i} and the references
given there and exercise 7.112 b in \cite{st}).

The number of planar trees with $n+1$ nodes and the number of planar
binary trees with $n+1$ leaves are both given by the $n$'th Catalan
number. It follows that the number of elements in $\mathcal{B}_n$ is
given by the above formula. 

\begin{cor}
The number $d(n)$ of elements in the mutation class of any quiver of
type $D_n$ is given by:
\[
d(n) = \left\{ \begin{array}{l l}
  \sum_{d|n} \phi(n/d)\binom{2d}{d}/(2n) & \quad \mbox{if $n \geq 5$}, \\
  6 & \quad \mbox{if $n=4$},
\end{array} \right. 
\]
where $\phi$ is the Euler function. 
\end{cor}

We proved this for $n \geq 5$ and for $n=4$ the number is $6$. See Figure \ref{figmutationclassD4} for all quivers in the mutation class of $D_4$. See Table \ref{tableexamples} for some values of $d(n)$.

\begin{table}
\begin{tabular}{l|l}
$n$&$d(n)$\\ \hline
3&4\\
4&6\\
5&26\\
6&80\\
7&246\\
\end{tabular}
\; \;
\begin{tabular}{l|l}
$n$&$d(n)$\\ \hline
8&810\\
9&2704\\
10&9252\\
11&32066\\
12&112720\\
\end{tabular}
\caption{Some values of $d(n)$.}\label{tableexamples}
\end{table}

\section{Mutation of trees}
We want to define a mutation operation on the elements in $\mathcal{B}$,
and we want this to commute with mutation of triangulations. Mutating a triangulation at a given diagonal is defined as removing this diagonal and replacing it with another one to obtain a new triangulation. This can be done in one and only one way. 

Let $\Delta$ be a triangulation in $\mathcal{T}_n$ and let $\sigma(\Delta) = T$ be the corresponding tree. An inner edge of $T$ corresponds to a diagonal in $\Delta$ not going from the puncture to the border, since the edges crosses these diagonals when we construct $T$ from $\Delta$. However, when we construct $T$ from $\Delta$, no edges in $T$ crosses a diagonal between the puncture and the border. To define mutation on $T$ corresponding to mutating at a diagonal $\alpha$ between the puncture and the border, we instead define mutation at two adjacent edges from the root in $T$, namely the two edges from the root in $T$ separated by $\alpha$. 

\begin{figure}[htp]
  \begin{center}
    \includegraphics[width=12.5cm]{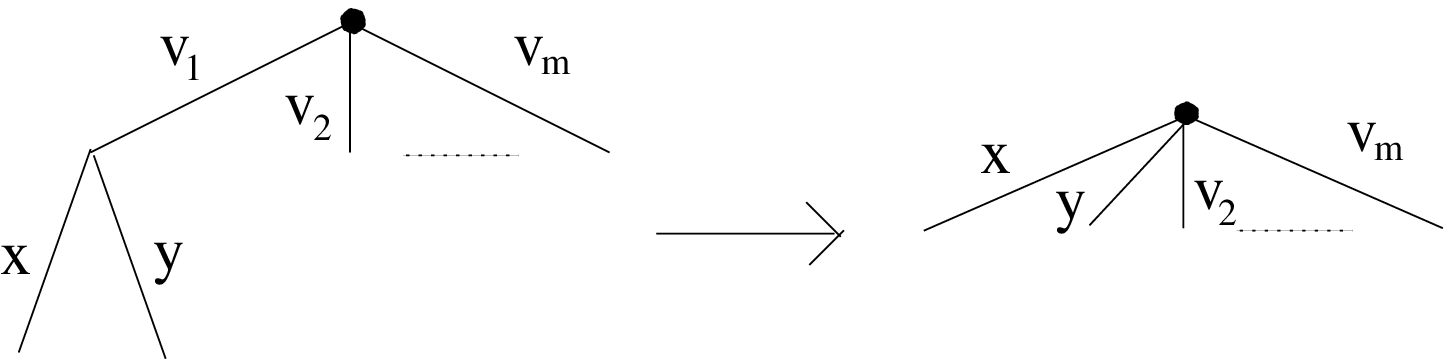}
    \includegraphics[width=12.5cm]{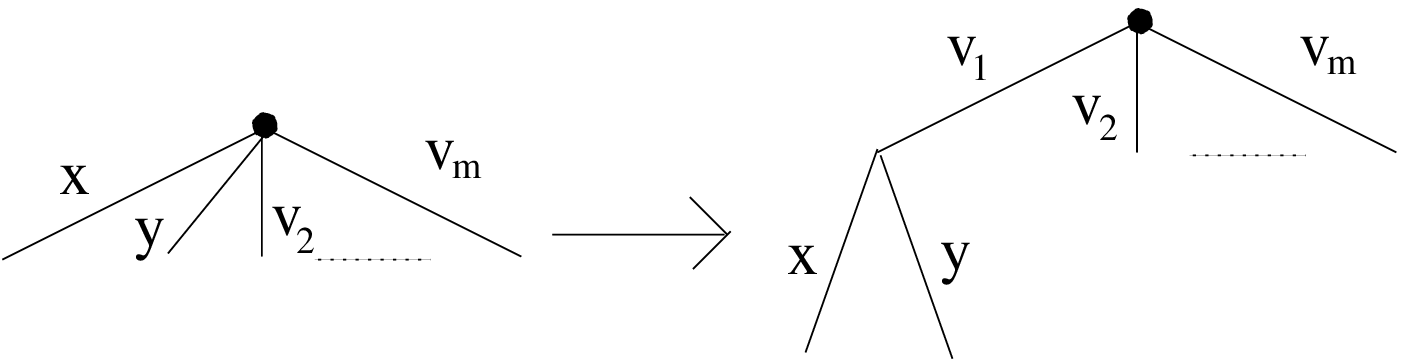}
    \includegraphics[width=12.5cm]{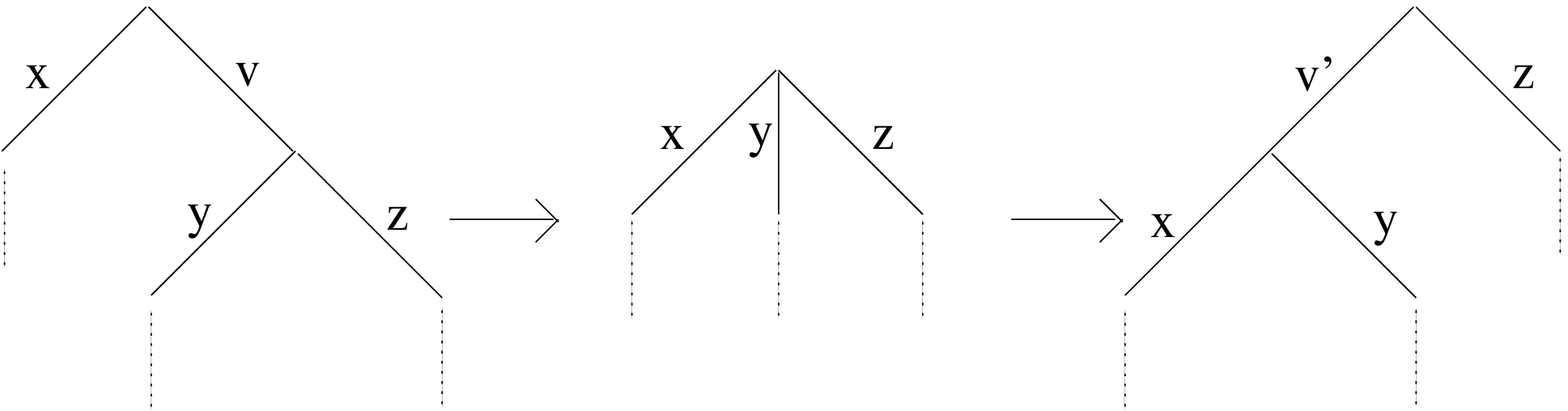}
    \includegraphics[width=12.5cm]{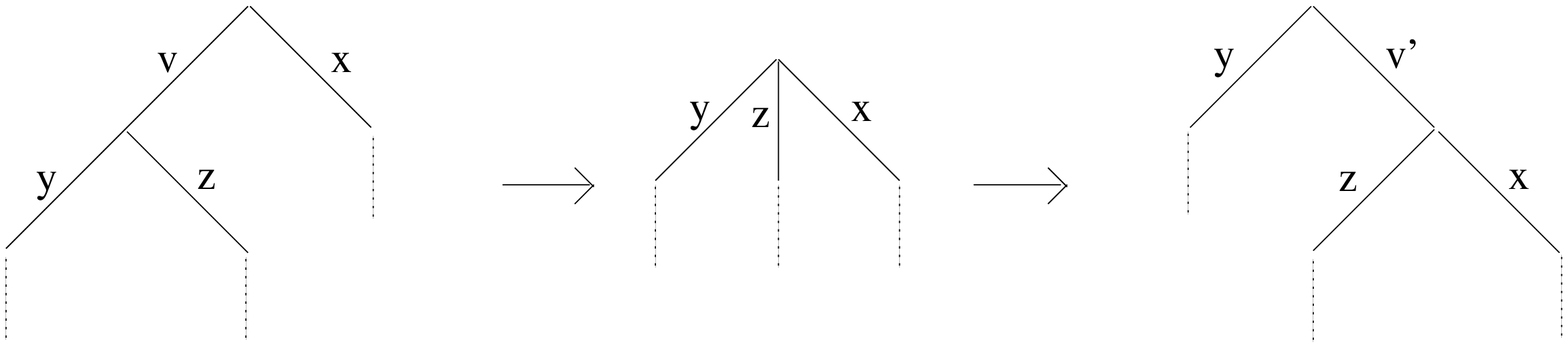}
  \end{center}\caption{Mutation of a tree.}
  \label{figtreemutate1}
  \end{figure}

\begin{enumerate}   
\item Let $v_1$ be an edge from the root in $T$. The mutation of $T$ at
  $v_1$ is a new tree obtained in the following way. Remove the edge
  $v_1$. Identify the root of the full subtree of $T$ ending in $v_1$ with the
  root in $T$. See the first picture in Figure \ref{figtreemutate1}.
\item Let $x$ and $y$ be two adjacent edges from the root of
  $T$. The mutation of $T$ at $x$ and $y$ is a new tree obtained in
  the following way. Disconnect the full subtree of $T$ containing $x$
  and $y$. Add an edge $v_1$ from the root and connect the subtree to
  the end of $v_1$. See the second picture in Figure \ref{figtreemutate1}. 
\item Let $v$ be an inner edge not going from the root or to a
  leaf. The mutation of $T$ at $v$ is a new tree obtained in the
  following way. Suppose $v$ is an edge from the nodes $r$ to $t$, going down in
  the tree. Let $x$ be the other edge starting in $r$, and let $y$ and
  $z$ be the two edges starting in $t$. See the third and fourth
  picture in Figure \ref{figtreemutate1}. Suppose $x$ goes
  to the left from $r$ and $v$ goes to the right, as in the third
  picture. Disconnect the full subtree with $t$ as a root. Remove the
  edge $v$ and identify $r$ with $t$. Disconnect the full subtree $T'$
  containing $x$ and $y$. Create a new vertex $v'$ starting in $r$ and
  identify the root of $T'$ with the node ending in $v'$. See the
  third picture in Figure \ref{figtreemutate1}. If $x$ goes to the
  right from $r$ and $v$ goes to the left, we define mutation at $v$
  in a similar way as shown in the fourth picture.  
\end{enumerate}

We claim that mutation of a tree as defined above commutes with
mutation of triangulations. We leave the details of the proof to the reader. 

\begin{prop}\label{treetricommutes}
Mutation of trees commutes with mutations of triangulations and
quivers. 
\end{prop}
\begin{proof}[Sketch of proof:]

For mutation of type 3, we mutate at a diagonal not going between the puncture and the border, so we are in the situation shown in Figure
\ref{muttype3}. We see that mutation of trees commutes with
mutation of triangulations. 

\begin{figure}[htp]
  \begin{center}
    \includegraphics[width=10cm]{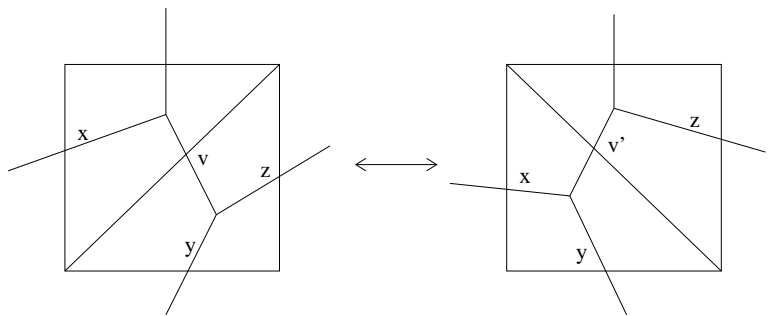}
  \end{center}\caption{Mutation of triangulations and trees commute. See proof of Proposition \ref{treetricommutes}.}
  \label{muttype3}
  \end{figure}

For mutation of type 1 and 2 we have to consider the three cases shown in Figure \ref{muttype1}. In these cases we also see that mutation as defined above commutes with mutation of triangulations.   

\begin{figure}[htp]
  
	\begin{minipage}[b]{12.5 cm}
	\hspace{1.5 cm}
    \includegraphics[width=10cm]{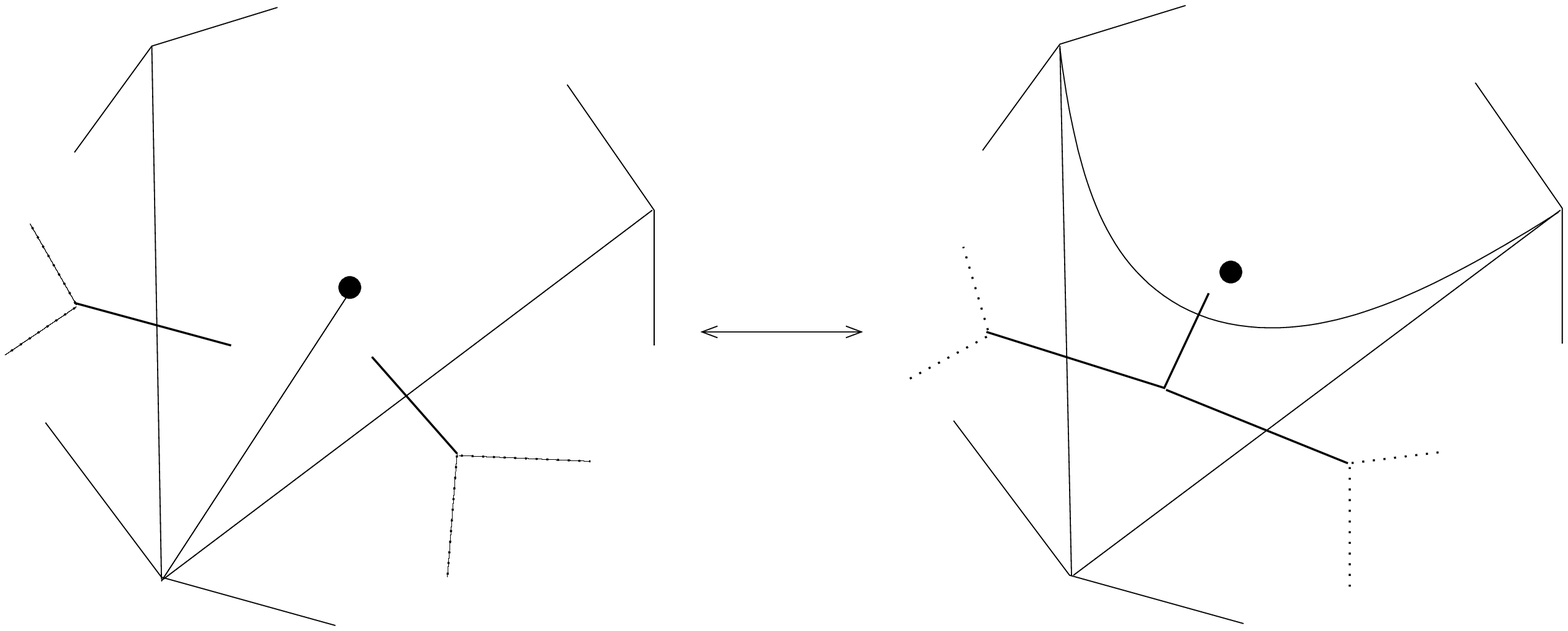} 
\vspace{1 cm}
	\end{minipage}

    	\begin{minipage}[b]{12.5 cm}
\hspace{1.5 cm}    
\includegraphics[width=10cm]{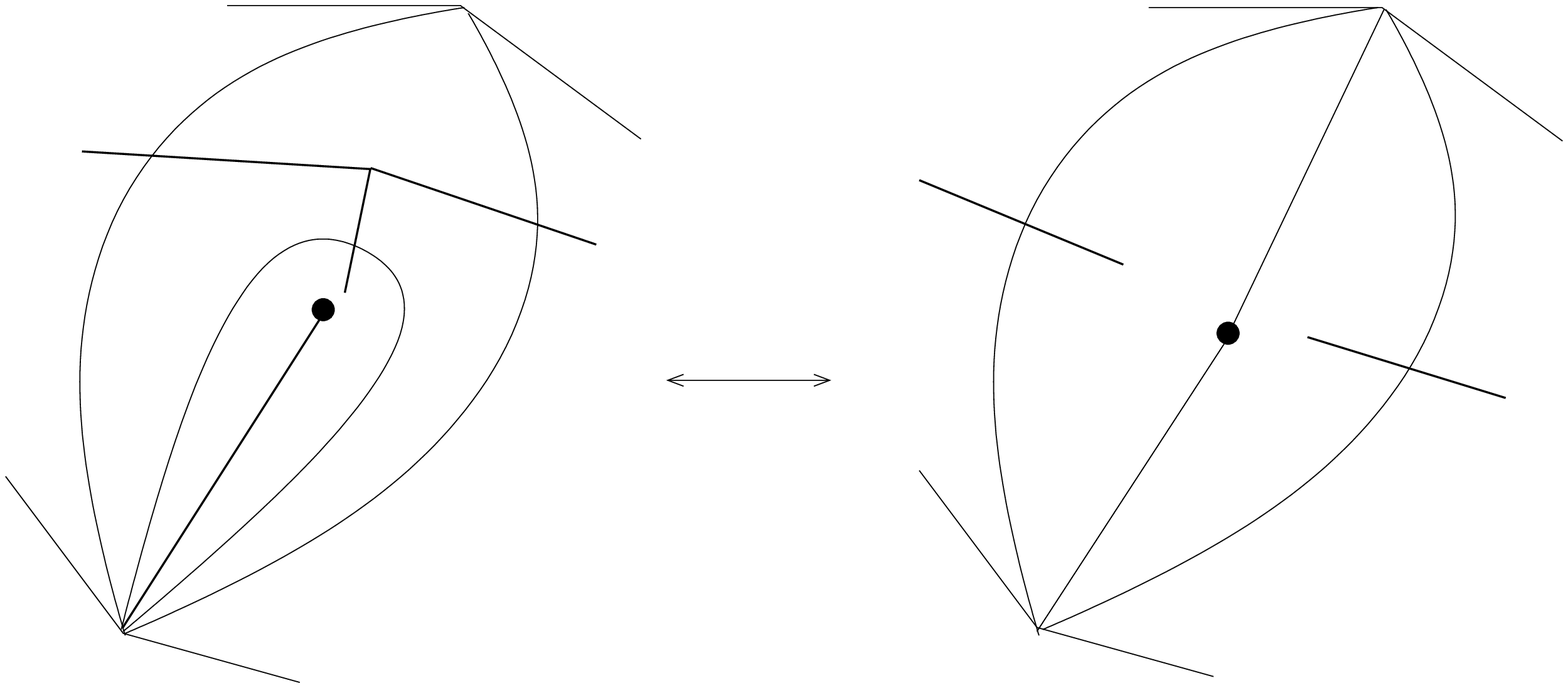}
\vspace{1 cm}
    	\end{minipage}

	\begin{minipage}[b]{12.5 cm}
\hspace{1.5 cm}    
\includegraphics[width=10cm]{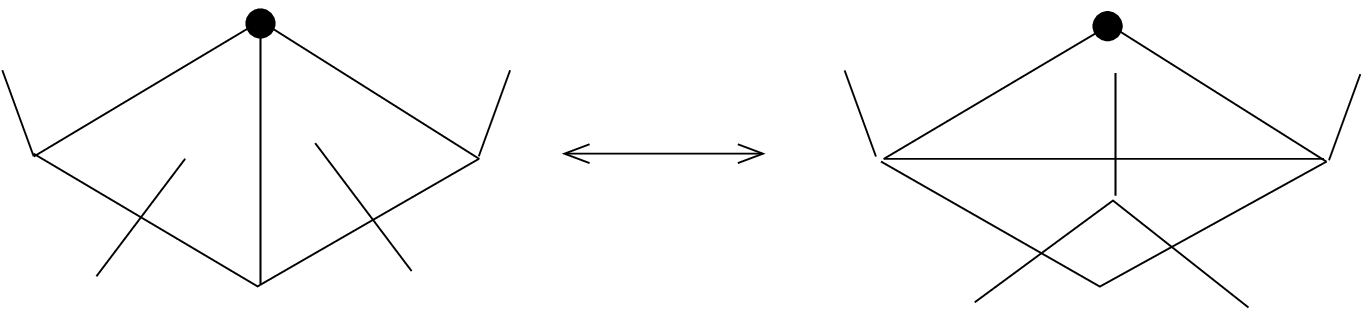}
	\end{minipage}
  \caption{Mutation of triangulations and trees commutes. See proof of Proposition \ref{treetricommutes}.}
  \label{muttype1}
  \end{figure}

\end{proof}

\begin{figure}\centering
\includegraphics[width=12.0cm]{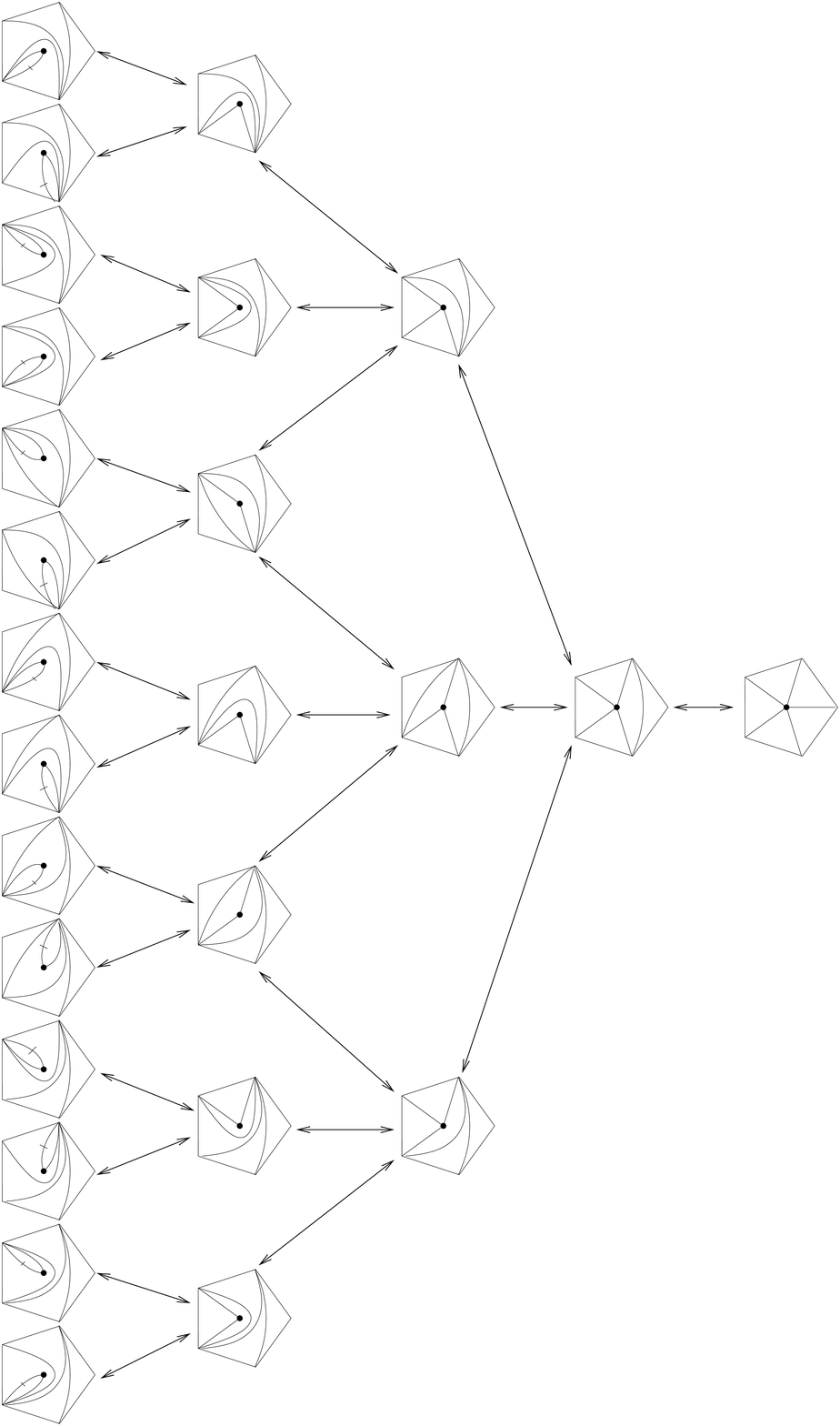}
\caption{All triangulations of type $D_5$.}
\label{triangulations tree}
\end{figure}

\begin{figure}\centering
\includegraphics[width=11.5cm]{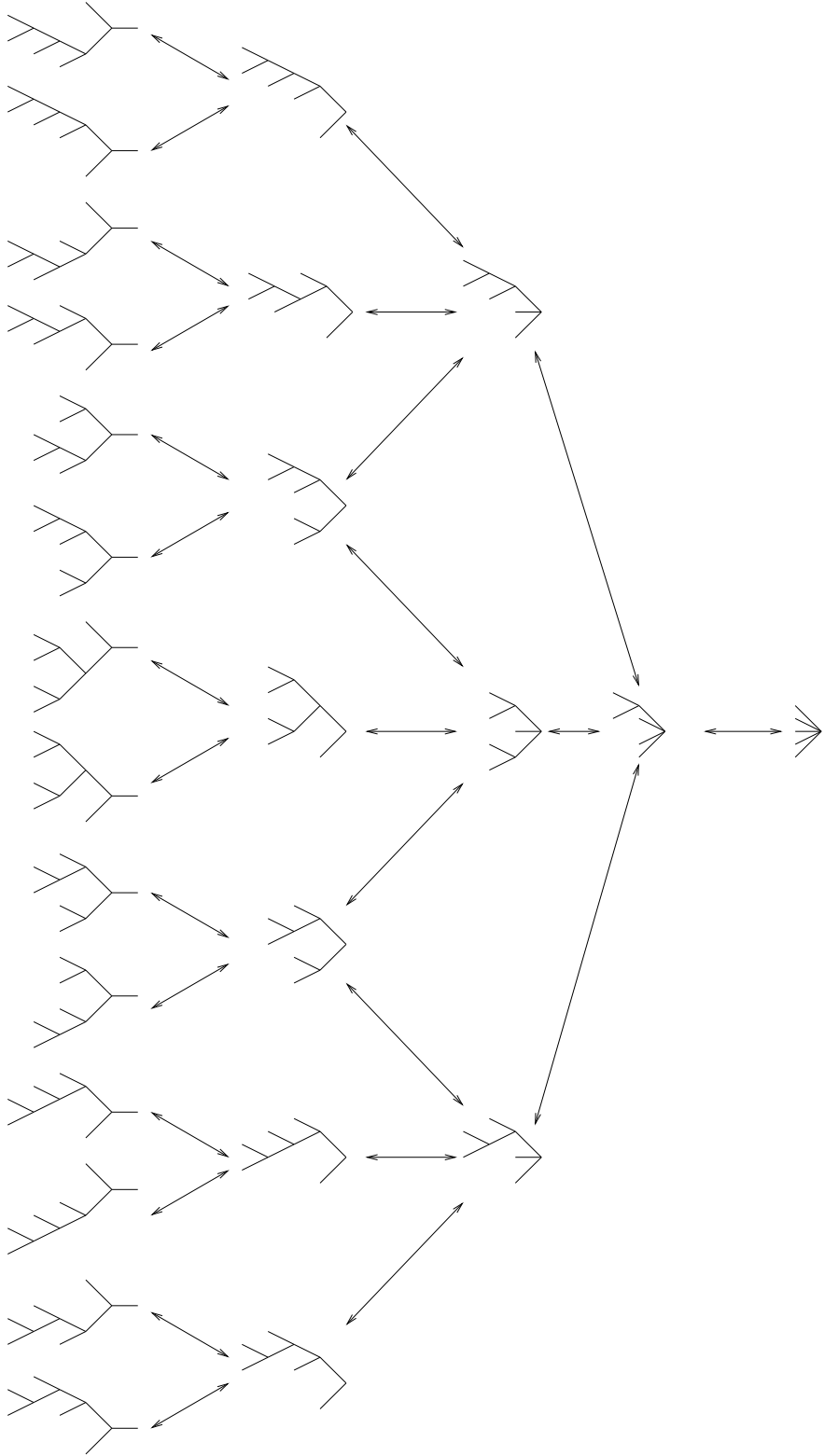}
\caption{All trees of type $D_5$.}
\label{tree tree}
\end{figure}

Figure \ref{triangulations tree} and \ref{tree tree} shows the
mutations of type 1 and 2 for both triangulations and trees in the
$D_5$ case. Note that mutation of type 2 adds an edge from the root, or equivalently replaces a diagonal not between the puncture and the border with a diagonal between the puncture and the border. Mutation of type 1 is the opposite operation. This defines a tree of mutations as shown in Figure \ref{triangulations tree} and \ref{tree tree}, where going down in the tree corresponds to mutatation of type 1 and going up in the tree corresponds to mutation of type 2. If we drew arrows for mutations of type 3, the arrows would go to trees (or triangulations) in the same level in the tree of mutations. It is easy to see that this holds in general for any $n$.  

\small

\end{document}